\documentclass[11pt]{article}
\usepackage{amssymb,amsmath,verbatim,amsthm}
\newtheorem{Theorem}{Theorem}[section]
\newtheorem{Lemma}[Theorem]{Lemma}
\newtheorem{Corollary}[Theorem]{Corollary}

\newtheorem{Construction}[Theorem]{Construction}

\newtheoremstyle{noparens}
  {}{}
  {\itshape}{}
  {\bfseries}{.}
  { }
  {\thmname{#1}\thmnumber{ #2}\mdseries\thmnote{ #3}}
\theoremstyle{noparens}
\usepackage{enumerate}
\usepackage{color,xcolor}
\usepackage{indentfirst}
\newcommand{\ignore}[1]{}

\begin{document}

\renewcommand{\theequation}{\thesection.\arabic{equation}}

\title{\bf Steiner Quadruple Systems with Minimum Colorable Derived Designs: Constructions and Applications\footnote{The research was supported
 by National Natural Science Foundation of China (12171028, 12371326).}}
\author{
{\small Yuli  Tan,}  {\small  Junling  Zhou}\\
{\small School of Mathematics and Statistics}\\
{\small Beijing Jiaotong University}\\
{\small Beijing  100044, China}\\
{\small YL.Tan@bjtu.edu.cn}\\
{\small jlzhou@bjtu.edu.cn}\\
}
\date{ }
\maketitle

\begin{abstract}
An $r$-block-coloring, simply $r$-coloring, of a Steiner triple system $\mathrm{STS}(v)$ is a partition of the block set into $r$ color classes, each color class being a partial parallel class.
The chromatic index of $\mathrm{STS}(v)$, denoted by $\chi^{\prime}(v)$, is the smallest $r$ for which an $r$-coloring of an $\mathrm{STS}(v)$ exists.
A minimum colorable Steiner triple system $\mathrm{mcSTS}(v)$ is an $\mathrm{STS}(v)$ admitting a $\chi^{\prime} (v)$-coloring.
We generalize the notion of an $\mathrm{RDSQS}$ (a Steiner quadruple system $\mathrm{SQS}$ with resolvable derived designs) to $\mathrm{mcDSQS}$, 
representing an $\mathrm{SQS}$ whose derived design at every point is minimum colorable.
This is motivated from an application in non-binary diameter perfect codes.
The purpose of this paper is to display a few recursive constructions to produce $\mathrm{mcDSQS}$s via Steiner systems $\mathrm{S}(3,K,v)$ with certain properties.
Among others, a  construction for $\mathrm{mcDSQS}$s is developed, which is also new even for $\mathrm{RDSQS}$s; special constructions concentrating only on $\mathrm{mcDSQS}(6n+2)$s are  demonstrated as well. As the main results, both
a new infinite family of $\mathrm{RDSQS}(6n+4)$s and the first infinite family of $\mathrm{mcDSQS}(6n+2)$s are constructed. To be specific, an $\mathrm{RDSQS}(2^{2m+1}+2)$
and an $\mathrm{mcDSQS}(2\cdot 9^{m}+2)$ are proved to exist, in which the former class gives rise to a new infinite family of large sets of Kirkman triple systems.
As applications, the smallest $q$ is  determined  such that a diameter perfect constant-weight $(n,\frac{1}{4}\tbinom{n}{3},6;4)_{q}$ code exists where
 $n \in\{ 2\cdot 9^{m}+2:m\geq 1\}\bigcup\{ 2^{2m+1}+2:m\geq 0\}$.
\medskip

\medskip
\end{abstract}

\section{Introduction and preliminaries}

Let $X$ be a set of $v$ points and $K$ a set of positive integers.
Let $\mathcal{G}$ be a collection of subsets (called {\em{groups}} or {\em{holes}}) of $X$ which partition $X$.
A {\em{group divisible $t$-design}} of order $v$ with block sizes from $K$, denoted by $\mathrm{GDD}(t,K,v)$, is a triple $(X,\mathcal{G},\mathcal{B})$ with the following properties:\vspace{-0.2cm}
\begin{itemize}
  \item [(1)] $\mathcal{B}$ is a family of subsets of $X$ (called {\em{blocks}}), each having cardinality in $K$;\vspace{-0.2cm}
 \item [(2)] each block of $\mathcal{B}$ intersects any given group in at most one point;\vspace{-0.2cm}
 \item [(3)] each $t$-subset of $X$ from $t$ distinct groups is contained in exactly one block.
\end{itemize}
\vspace{-0.2cm}
The {\em{type}} of the $\mathrm{GDD}$ is defined to be the multiset $T=\{|G|: G \in \mathcal{G}\}$ of group sizes.
We also use $g_{1}^{a_{1}} g_{2}^{a_{2}}\ldots g_{r}^{a_{r}}$ to denote the type, which means that there are $a_{i}$ groups of size $g_{i}$, $1\leq i \leq r$.
A $\mathrm{GDD}(t,\{k\},v)$ is also denoted by $\mathrm{GDD}(t,k,v)$.
A $\mathrm{GDD}(t,k,kn)$ of type $n^{k}$ is known as a {\em{transversal $t$-design}} and denoted by $\mathrm{TD}(t,k,n)$.
If $(X,\mathcal{G},\mathcal{B})$ is a $\mathrm{GDD}(t,K,v)$ of type $1^{v}$, then we omit the notation of $\mathcal{G}$ and call $(X,\mathcal{B})$ a {\em{Steiner $t$-system}}, also denoted $\mathrm{S}(t,K,v)$. The notation $\mathrm{S}(t,k,v)$ is used for $K=\{k\}$.
An $\mathrm{S}(2,3,v)$ is a {\em{Steiner triple system}} $\mathrm{STS}(v)$ and an $\mathrm{S}(3,4,v)$ is a {\em{Steiner quadruple system}} $\mathrm{SQS}(v)$.
It is well-known that an $\mathrm{STS}(v)$ exists if and only if $v \equiv 1,3 \pmod{6}$ \cite{TSTS} and an $\mathrm{SQS}(v)$ exists if and only if $v \equiv 2,4 \pmod{6}$ \cite{HSQS}.

Let $(X,\mathcal{G},\mathcal{B})$ be a $\mathrm{GDD}$ and $Y\subseteq X$.
A {\em{partial parallel class}} ($\mathrm{PPC}$) of $Y$ is a collection of pairwise disjoint blocks contained in $Y$.
A {\em{parallel class}} ($\mathrm{PC}$) of $Y$ is a $\mathrm{PPC}$ that partitions $Y$.
We say that $(X, \mathcal{G}, \mathcal{B})$ is {\em {resolvable}}, denoted by $\mathrm{RGDD}$, if $\mathcal{B}$ can be partitioned into $\mathrm{PC}$s of $X$.
An $\mathrm{RGDD}(t,k,v)$ of type $1^{v}$ is a resolvable $\mathrm{S}(t,k,v)$, denoted by $\mathrm{RS}(t,k,v)$. An $\mathrm{RS}(2,3,v)$ is a {\em{Kirkman triple system}} $\mathrm{KTS}(v)$.
A $\mathrm{KTS}(v)$ exists if and only if $v\equiv  3 \pmod{6}$ \cite{KTS1,KTS2}. Two $\mathrm{KTS}(v)$s on the same set $X$ are said to be {\em{disjoint}} if they have no blocks in common. A {\em{large set of Kirkman triple system}} of order $v$, denoted by $\mathrm{LKTS}(v)$, is a partition of all triples of a $v$-element set $X$ into $v-2$ pairwise disjoint $\mathrm{KTS}(v)$s.
Although the problem of the existence of $\mathrm{LKTS}(v)$s has a long history, it remains pretty much open.

Let $(X,\mathcal{B})$ be an $\mathrm{S}(t,k,v)$. For any point $x\in X$, let $\mathcal{B}_{x} = \{B\setminus \{x\}:x\in B, B\in \mathcal{B}\}$. Then $(X\setminus\{x\},\mathcal{B}_{x})$ is an $\mathrm{S}(t-1,k-1,v-1)$, which is said to be the {\em{derived design}} of $(X,\mathcal{B})$ at the point $x$.
An $\mathrm{S}(t,k,v)$ with {\em resolvable derived designs}, abbreviated to $\mathrm{RDS}(t,k,v)$, refers to an $\mathrm{S}(t,k,v)$ whose derived design at every point is resolvable.
We use the notation $\mathrm{RDSQS}(v)$ to denote an $\mathrm{SQS}(v)$ with resolvable derived designs. $\mathrm{RDSQS}$s play an important role in recursive constructions of large sets of Kirkman triple systems (see \cite{Zhou-lkts,LEI1,Kang-2008}).

\begin{Theorem}[\label{RDS_1}{\cite{Zhou-lkts,Kang-2008}}]
If there exists an $\mathrm{RDSQS}(v + 1)$, then there exists an $\mathrm{LKTS}(3v)$.
\end{Theorem}

We know that the necessary condition for the existence of an $\mathrm{RDSQS}(v)$ is $v \equiv 4 \pmod{6}$ as a $\mathrm{KTS}(v)$ exists if and only if $v\equiv  3 \pmod{6}$.
 An $\mathrm{SQS}(v)$ with $v\equiv 2\pmod 6$ dose not have a resolvable derived design, but one may instead ask for a derived design with a partition of the blocks into minimum number of PPCs.
We introduce the concept of block-coloring.
An {\em{$r$-block-coloring}}, or simply {\em{$r$-coloring}}, of a design is a partition of the block set into $r$ color classes, each color class being a $\mathrm{PPC}$.
The collection of the $r$ $\mathrm{PPC}$s is called a {\em{resolution}}. The {\em{chromatic index}} of a design $\mathcal{D}$, denoted by $\chi^{\prime}(\mathcal{D})$, is the smallest $r$ for which $\mathcal{D}$ admits an $r$-coloring. The chromatic index of Steiner systems $\mathrm{S}(t,k,v)$, denoted by $\chi^{\prime}_{t,k}(v)$, is defined as
$$\chi^{\prime}_{t,k}(v) = \mathop{\min}_{\mathcal{D}\text{ is an }S(t,k,v)}\chi^{\prime}(\mathcal{D}).$$
For convenience, the notation $\chi^{\prime} (v)$ is used for $\chi^{\prime}_{2,3}(v)$.
If $\mathcal{D}$ is a resolvable $\mathrm{S}(t,k,v)$, then $\chi^{\prime}_{t,k}(\mathcal{D})$ equals  ${k\tbinom{v}{t}}/{v\tbinom{k}{t}}$. Therefore, $\chi^{\prime}(v)$ equals $(v-1)/2$ for $v\equiv  3 \pmod{6}$ because of the existence of a $\mathrm{KTS}(v)$.
The unique $\mathrm{STS}(7)$ has chromatic index $7$ as each pair of blocks has an intersection point.
There are two nonisomorphic $\mathrm{STS}(13)$s and each of them has chromatic index $8$ (see \cite{HTS14}). Generally, $\chi^{\prime}(v)$ equals $(v+1)/2$ if $v\equiv  1 \pmod{6}$ and $v\notin \{7,13\}$ (see \cite{HTS}). To sum up,
\begin{align}\label{E1}
\chi^{\prime}(v) =
\begin{cases}
\frac{v-1}{2}, & \text{ if } v\equiv  3 \pmod{6}, \\
\frac{v+1}{2}, & \text{ if } v\equiv  1 \pmod{6}\text{ and } v\notin \{7,13\},\\
7, & \text{ if } v=7,\\
8, & \text{ if } v=13.\\
\end{cases}
\end{align}

If an $\mathrm{STS}(v)$ admits a $\chi^{\prime} (v)$-coloring, then it is said to be
{\em{minimum colorable}}, denoted by $\mathrm{mcSTS}(v)$.
If $v\equiv  3 \pmod{6}$, then an $\mathrm{mcSTS}(v)$ is a $\mathrm{KTS}(v)$.
An $\mathrm{SQS}(v)$ with {\em{minimum colorable derived designs}}, abbreviated to $\mathrm{mcDSQS}(v)$, refers to an $\mathrm{SQS}(v)$ whose derived design at every point is minimum colorable. Therefore, an $\mathrm{mcDSQS}(v)$ is the same as an $\mathrm{RDSQS}(v)$ if $v\equiv  4 \pmod{6}$. There is a large literature studying $\mathrm{RDSQS}$s \cite{zheng,Zhou-lkts,LEI1,yuanotripling}, while we only have one known $\mathrm{mcDSQS}(6n+2)$ with $n\geq 2$, that is, an $\mathrm{mcDSQS}(20)$ \cite{SHI}. In this work we introduce the concept of $\mathrm{mcDSQS}$s, not only as a generalization of $\mathrm{RDSQS}$s, but also due to their applications in diameter perfect codes.

Diameter perfect codes are a natural generalization of perfect codes. They are based on the code-anticode bound, a generalization of sphere-packing bound. Code-anticode bound was first introduced by Delsarte (see \cite{Delsarte}) and generalized by Ahlswede et al. (see \cite{DPD}).
Let $\mathcal{C}$ be a code in a metric space $\mathcal{V}$ with a transitive automorphism group, where the minimum distance between the codewords in $\mathcal{C}$ is $d$. Let $\mathcal{A}$ be a anticode of $\mathcal{V}$ such that the maximum distance (diameter) between elements in $\mathcal{A}$ is $D$, where $D< d$. Then
$$|\mathcal{C}|\cdot |\mathcal{A}|\leq |\mathcal{V}|.$$
If the inequality holds with equality, then $\mathcal{C}$ is said to be a {\em{$D$-diameter perfect code}} or {\em{diameter perfect distance-$d$ code}}.
Let $\mathcal{J}_{q}(n,w)$ denote the set of
all words of length $n$ and weight $w$, over $\mathbb{Z}_{q}=\{0,1,\ldots, q-1\}$, where the Hamming distance is used as the metric.
A code $\mathcal{C}\subseteq \mathcal{J}_{q}(n,w)$ is a {\em{$q$-ary constant-weight-$w$ code of length $n$}} and it is denoted by $(n,|\mathcal{C}|,d;w)_{q}$ code if $\mathcal{C}$ has minimum distance $d$.

In \cite{Etzion}, Etzion classified the known non-binary diameter perfect constant-weight codes into six families. Shi et al. \cite{SHI} generalized one of these families and investigated a family of such codes, where the supports of the codewords from a Steiner system $\mathrm{S}(t,w,n)$.
For $t\leq w$ and $2w-t\leq n$, define
$$\mathcal{A}_{q}(n,w,t)=\{(a_{1},\ldots,a_{t},\ldots,a_{n})\in \mathbb{Z}_{q}^{n}:\text{wt}((a_{1},\ldots,a_{t}))=t,\text{wt}((a_{t+1},\ldots,a_{n}))=w-t\}.$$
We can check that $\mathcal{A}_{q}(n,w,t)$ is a diameter-$(2w-t)$ anticode with ${{n-t}\choose{w-t}}(q-1)^{w}$ words in $\mathcal{J}_{q}(n,w)$.

\begin{Theorem}[\label{shi1}{\cite[Theorem 3]{SHI}}]
Let $t,w$ and $n$ be integers such that $0 <t \leq w < n$, $2w-t\leq n$.
\begin{itemize}
  \item [$(1)$] An $(n,\tbinom{n}{t}/\tbinom{w}{t},2w-t+1;w)_{q}$ code is a diameter perfect constant-weight code which attains the code-anticode bound with respect to the anticode $\mathcal{A}_{q}(n,w,t)$;\vspace{-0.2cm}
  \item [$(2)$] If there exists a Steiner system $\mathrm{S}(t,w,n)$, then there exists an $(n,\tbinom{n}{t}/\tbinom{w}{t},2w-t+1;w)_{q}$ code with $q\geq \tbinom{n-1}{t-1}/\tbinom{w-1}{t-1}+1$.
\end{itemize}
\end{Theorem}
Define $q_{0}^{\prime}(t,w,n)$ to be the smallest $q$ such that an $(n,\tbinom{n}{t}/\tbinom{w}{t},2w-t+1;w)_{q}$ code exists. Theorem \ref{shi1} (2) provides an upper bound of $q_{0}^{\prime}(t,w,n)$. When $t=3$ and $w=4$, the exact value of $q_{0}^{\prime}(t,w,n)$ is presented in terms of chromatic indexes of derived designs of an $\mathrm{SQS}(n)$.

\begin{Theorem}[\label{shi2}{\cite[Theorem 5]{SHI}}]
For $n\equiv 2,4 \pmod6$,
$$q_{0}^{\prime}(3,4,n)= \mathop{\min}_{S \text{ is } \mathrm{SQS}(n)}\mathop{\max}_{D\in S^{\prime}}\chi^{\prime}(D)+1,$$
where $S^{\prime}$ is the set of all $\mathrm{STS}(n-1)$s derived from $S$, the $\mathrm{SQS}(n)$.
\end{Theorem}
The following corollary is immediate.
\begin{Corollary}\label{shi-c1}
\begin{itemize}
 \item [$(1)$] For $n\equiv 2 \pmod6$ and $n\geq 20$, if there exists an $\mathrm{mcDSQS}(n)$, then $q_{0}^{\prime}(3,4,n)=\frac{n}{2}+1$.
 \item [$(2)$] For $n\equiv 4 \pmod6$, if there exists an $\mathrm{RDSQS}(n)$, then $q_{0}^{\prime}(3,4,n)=\frac{n}{2}$.
 \end{itemize}
\end{Corollary}

The rest of the paper is organized as follows. Section 2 is devoted to demonstrating recursive constructions for generic $\mathrm{mcDSQS}$s. In Subsection 2.1, we define a $\mathrm{gcDCQS}(g^{n}:s)$ to build an $\mathrm{mcDSQS}(gn+s)$ by filling in holes. Then we demonstrate several recursive constructions to produce $\mathrm{gcDCQS}(g^{n}:s)$ via Steiner systems $\mathrm{S}(3,K,v)$ with certain properties.
In Subsection 2.2, we concentrate on constructions of $\mathrm{mcDSQS}(6n+2)$s. We propose two new types of  colorings of the derived designs of a $\mathrm{CQS}(g^{n}:2)$ and then present related constructions.
In Section 3, we construct some necessary designs required in our recursive constructions. Then we establish the existence of a new infinite family of $\mathrm{RDSQS}(6n+4)$s and  an infinite family of $\mathrm{mcDSQS}(6n+2)$s. 
As  consequences, we show their applications in diameter perfect codes and large sets of Kirkman triple systems. In Section 4, we conclude the paper with some open questions.

\section{Recursive constructions}

A combinatorial structure named candelabra quadruple system $\mathrm{CQS}$ plays a central role in the recursive constructions of $\mathrm{SQS}$s. In order to produce
$\mathrm{RDSQS}$s, candelabra quadruple systems with resolvable derived designs $\mathrm{RDCQS}$s were introduced and intensively studied in the last decade. 
In this section, we introduce several coloring types of the derived designs of a $\mathrm{CQS}$ as a generalization of an $\mathrm{RDCQS}$. Then we display a number of recursive constructions to produce $\mathrm{mcDSQS}$s. For later use, we introduce or recall some definitions.

An {\em{incomplete Steiner triple system}} on a $v$-set $X$ with an $h$-subset $H$ as a hole, denoted by $\mathrm{STS}(v,h)$, is a triple $(X,H,\mathcal{B})$, where $\mathcal{B}$ is a block set of a $\mathrm{GDD}(2,3,v)$ of type $1^{v-h}h^{1}$ with $H$ as the long group.
Further let $v\equiv 1,3\pmod 6$.
If $\mathcal{B}$ can be partitioned into $r:=\chi^{\prime}(v)$ $\mathrm{PPC}$s of $X$, in which there are $s:=\chi^{\prime}(v)-\frac{v-h}{2}$ $\mathrm{PPC}$s of $X\setminus H$,
then we say that $(X,H,\mathcal{B})$ is {\em{good colorable}}, denoted by $\mathrm{gcSTS}(v,h)$.
It then follows from (\ref{E1}) that
\begin{align*}
(r,s) = \begin{cases}
(\frac{v-1}{2},\frac{h-1}{2}), & \text{ if } v\equiv  3 \pmod{6},\\
(\frac{v+1}{2},\frac{h+1}{2}), & \text{ if } v\equiv 1 \pmod{6}, v\neq 7,13.\\
		   \end{cases}
\end{align*}
It is not difficult to check that, when $v\equiv 3\pmod 6$, a $\mathrm{gcSTS}(v,h)$ $(X,H,\mathcal{B})$ becomes an {\em incomplete Kirkman triple system }$\mathrm{KTS}(v,h)$, whose blocks are partitioned into $\frac{v-h}{2}$ $\mathrm{PC}$s of $X$ and $\frac{h-1}{2}$ $\mathrm{PC}$s of $X\setminus H$.

An {\em{incomplete  Steiner quadruple system}} on a $v$-set $X$ with an $h$-subset $H$ as a hole, denoted by $\mathrm{SQS}(v,h)$, is a triple $(X,H,\mathcal{B})$, where $\mathcal{B}$ is a set of quadruples of $X$ (blocks), such that every triple $T\subseteq X$ with $T\not\subseteq H$ is contained in a unique block and no triple of $H$ is contained in any block.
We can also define the derived designs of an $\mathrm{SQS}(v,h)$.
For $x \in X$, put $\mathcal{B}_x=\{B \setminus\{x\}: x \in B, B \in \mathcal{B}\}$.
It is immediate that each $(X \setminus\{x\},\mathcal{B}_x)$ with $x\in X\setminus H$ forms an $\mathrm{STS}(v-1)$ and each $(X \setminus \{x\}, H\setminus\{x\}, \mathcal{B}_x)$ with $x \in H$ forms an $\mathrm{STS}(v-1,h-1)$.
An $\mathrm{SQS}(v,h)$ with {good colorable derived designs}, denoted by $\mathrm{gcDSQS}(v,h)$, if all derived $\mathrm{STS}(v-1)$s are minimum colorable and all derived $\mathrm{STS}(v-1,h-1)$s are good colorable.

A {\em{Kirkman frame}} $\mathrm{F}(2,3,v\{m\})$ (as in \cite{frame}), briefly by $\mathrm{KF}(m^{v})$, is a $\mathrm{GDD}(2,3,mv)$ of type $m^{v}$ $(X,\mathcal{G},\mathcal{B})$ such that the block set $\mathcal{B}$ can be partitioned into subsets $P_{r}$, $r\in R$, each $P_{r}$ being a $\mathrm{PC}$ of $X\setminus G$ for some group $G\in \mathcal{G}$. A $\mathrm{KF}(m^{v})$ exists if and only if $m\equiv 0 \pmod 2$, $v\geq 4$ and $m(v-1)\equiv 0 \pmod 3$ (see \cite{KF}).

Let $(X, \mathcal{G}, \mathcal{A})$ be a $\mathrm{GDD}(3,4, mv)$ of type $m^v$. For any $x \in X$, let $\mathcal{A}_x=\{A \backslash\{x\}: x \in A, A \in \mathcal{A}\}$. Then the derived design $(X \backslash G, \mathcal{G} \backslash\{G\}, \mathcal{A}_x)$ is a $\mathrm{GDD}(2,3, m(v-1))$ of type $m^{v-1}$.
Here we examine two types of coloring of the derived designs. We use $\mathrm{RDGDD}(3,4, v\{m\})$ to denote a $\mathrm{GDD}(3,4, mv)$ of type $m^{v}$ with resolvable derived designs, meaning that its derived design at every point is a resolvable $\mathrm{GDD}(2,3,m(v-1))$.
An $\mathrm{RDTD}(3,4,m)$ denotes an $\mathrm{RDGDD}(3,4,4\{m\})$.
Parallelly, we use $\mathrm{FDGDD}(3,4,v\{m\})$ to denote a $\mathrm{GDD}(3,4,mv)$ of type $m^{v}$ whose derived design at every point is a $\mathrm{KF}(m^{v-1})$.

Let $X$ be a $v$-element set. Suppose that $\mathcal{G}$, $\mathcal{B}_{i}$ $(1\leq i \leq s)$ and $\mathcal{T}$ are collections of some subsets of $X$ and $\mathcal{G}$ is a partition of $X$.
Let $K_{i}$ $(1\leq i \leq s)$ and $K_{T}$ be the block size sets of $\mathcal{B}_{i}$ and $\mathcal{T}$ respectively.
An {\em{$s$-fan design}}, denoted by $s$-$\mathrm{FG}(3,(K_1, \ldots, K_s, K_T), v)$, is an $(s+3)$-tuple $(X,\mathcal{G},\mathcal{B}_{1},\ldots,\mathcal{B}_{s},\mathcal{T})$ such that each $(X,\mathcal{G}\cup \mathcal{B}_{i})$ with $1\leq i \leq s$ is a Steiner $2$-system and $(X,\mathcal{G}\bigcup (\bigcup_{i=1}^{s}\mathcal{B}_{i})\bigcup \mathcal{T})$ is a Steiner $3$-system. 
The type of an $s$-$\mathrm{FG}$ is the multiset $\{|G|: G \in \mathcal{G}\}$. Note that an $s$-$\mathrm{FG}(3,(K_1, \ldots, K_s, K_T), v)$
of type $1^v$ is an $\mathrm{S}(3,(\bigcup_{i=1}^{s}K_{i})\bigcup K_{T},v)$ having $s$ subdesigns $\mathrm{S}(2,K_{i},v)$. For such a case we simply say that $(X,\mathcal{B}_{1},\ldots,\mathcal{B}_{s},\mathcal{T})$ is an $s$-fan $\mathrm{S}(3,(K_{1},\ldots, K_{s},$ $K_{T}),v)$.

 For a $v$-set $X$, take an $s$-subset $S\subseteq X$ and call it a {\em{stem}}. Let $\mathcal{G}$ be a collection of subsets, called groups, of $X\setminus S$ which partition $X\setminus S$.
A {\em{candelabra $t$-system}} (as in \cite{Cande}) of order $v$ with block sizes from a set $K$, briefly by $\mathrm{CS}(t,K,v)$, is a quadruple $(X, S, \mathcal{G}, \mathcal{A})$ satisfying that
\vspace{-0.2cm}
\begin{itemize}
  \item [(1)] $\mathcal{A}$ is a family of subsets of $X$, called blocks, each of cardinality from $K$; \vspace{-0.2cm}
  \item [(2)] each $t$-subset  of $X$ except those contained in $S\cup G$ for some $G\in \mathcal{G}$ is contained in exactly one block (while no $t$-subset  in $S\cup G$ is contained in any block).
\end{itemize}
\vspace{-0.2cm}
By the type of a $\mathrm{CS}(t,K,v)$ we mean the list $(\{|G|: G \in \mathcal{G}\}: s)$. Exponential notation is also applied to denote the type of $\mathcal{G}$. A candelabra system with $t=3$ and $k=4$ is a {\em{candelabra quadruple system}} and denoted by $\mathrm{CQS}$. 

It is a standard construction to build an $\mathrm{SQS}(gn+s)$ from a $\mathrm{CQS}(g^{n}:s)$ by filling in holes. To produce $\mathrm{mcDSQS}$s, we need to define desirable
 colorings of the derived designs of  $\mathrm{CQS}$s. The rest of this section is divided into two subsections. In the first subsection, we generalize  known
  definitions and constructions related with $\mathrm{RDSQS}$s to develop recursive constructions applicable for generic $\mathrm{mcDSQS}$s. 
 In the second subsection, we turn our attention to $\mathrm{mcDSQS}$s only of orders $6n+2$. We introduce two new types of colorings of derived designs of a $\mathrm{CQS}(g^{n}:2)$ and display effective recursive constructions.

\subsection{Generic constructions}

Let $(X, S, \mathcal{G}, \mathcal{A})$ be a $\mathrm{CQS}(g^n: s)$.
For $x \in X$, put $\mathcal{A}_x=\{A \setminus\{x\}: x \in A, A \in \mathcal{A}\}$. It can be checked that each derived design $(X \setminus\{x\},(G \cup S) \setminus\{x\}, \mathcal{A}_x)$ with $x \in G$ and $G \in \mathcal{G}$ forms an $\mathrm{STS}(gn+s-1, g+s-1)$; and each derived design $(X \setminus S, \mathcal{G}, \mathcal{A}_x)$ with $x \in S$ forms a $\mathrm{GDD}(2,3, gn)$ of type $g^n$.
A $\mathrm{CQS}(g^n: s)$ is said to have {\em good colorable derived designs}, denoted by $\mathrm{gcDCQS}(g^n: s)$, if the collection $\{\mathcal{A}_x:x \in X\}$ of its derived designs has the following properties:\vspace{-0.2cm}
\begin{itemize}
  \item [(1)] for $x\in G$ and $G\in \mathcal{G}$, each $(X\setminus \{x\}, (G\cup S)\setminus \{x\},\mathcal{A}_{x})$ is a $\mathrm{gcSTS}(gn+s-1,g+s-1)$;\vspace{-0.2cm}
  \item [(2)] for $x\in S$, each $(X\setminus S, \mathcal{G}, \mathcal{A}_{x})$ is a $\mathrm{KF}(g^{n})$. \vspace{-0.2cm}
\end{itemize}

Whenever $gn+s\equiv 4\pmod 6$ and $s\ne 1$, a $\mathrm{gcDCQS}(g^{n}:s)$ becomes an $\mathrm{RDCQS}(g^{n}:s)$, the notation used in \cite{Zhou-lkts}.
A simple but important use of $\mathrm{gcDCQS}$s is the filling construction to build $\mathrm{mcDSQS}$s.
{\Construction \label{mcdcqs-mcdsqs}
If there exist a $\mathrm{gcDCQS}(g^n: s)$, an $\mathrm{mcDSQS}(g+s)$
 and a $\mathrm{gcDSQS}(g+s,s)$, then there exists an $\mathrm{mcDSQS}(gn+s)$.}

\proof 
See \cite{Zhou-lkts} for the case $gn+s \equiv 4 \pmod 6$.
We only need to handle the case $gn+s \equiv  2 \pmod 6$.
Suppose $(X,S,\mathcal{G},\mathcal{A})$ is the given $\mathrm{gcDCQS}(g^n: s)$.
For $x \in X$, let $\mathcal{A}_x=\{A \setminus\{x\}: x \in A, A \in \mathcal{A}\}$.
For $x\in G$ and $G\in \mathcal{G}$, each $(X\setminus \{x\}, (G\cup S)\setminus \{x\},\mathcal{A}_{x})$ is a $\mathrm{gcSTS}(gn+s-1,g+s-1)$; its resolution is $\{P_{x}^{j}:1 \leq j \leq \frac{gn+s}{2}\}$, where each
$P_{x}^{j}$ is a $\mathrm{PPC}$ of $X\setminus \{x\}$;
particularly, each $P_{x}^{j}$ with $1 \leq j \leq \frac {g+s}{2}$ is a $\mathrm{PPC}$ of $X\setminus (G\cup S)$.
For $x\in S$ , each $(X\setminus S, \mathcal{G}, \mathcal{A}_{x})$ is a $\mathrm{KF}(g^{n})$; its resolution is $\{P_{x}^{j}(G):G\in \mathcal{G}, 1\leq j\leq \frac{g}{2}\}$, where each $P_{x}^{j}(G)$ is a $\mathrm{PC}$ of $X\setminus (G\cup S)$.

Let $G_{0}\in \mathcal{G}$. We can construct by assumption an $\mathrm{mcDSQS}(g+s)$ on $G_{0}\cup S$ with block set $\mathcal{B}(G_{0})$. So we have $g+s$ derived $\mathrm{mcSTS}(g+s-1)$s, say $((G_{0}\cup S) \setminus\{x\}, \mathcal{B}_{x}(G_{0}))$ with $x\in G_{0}\cup S$.
Further, each $\mathcal{B}_{x}(G_{0})$ has a resolution $\{Q_{x}^{j}(G_{0}):1\leq j \leq \frac {g+s}{2}\}$ into $\mathrm{PPC}$s of $(G_{0}\cup S)\setminus \{x\}$.

For each $G\in \mathcal{G}$ with $G\neq G_{0}$, we construct a $\mathrm{gcDSQS}(g+s,s)$ on $G\cup S$ with $S$ as a hole and with block set $\mathcal{B}(G)$. So we have $g$ derived $\mathrm{mcSTS}(g+s-1)$s $((G\cup S) \setminus\{x\}, \mathcal{B}_{x}(G))$ with $x\in G$ and $s$ derived $\mathrm{gcSTS}(g+s-1,s-1)$s $((G\cup S) \setminus\{x\}, S\setminus \{x\},\mathcal{B}_{x}(G))$ with $x\in S$.
Further, each $\mathcal{B}_{x}(G)$ with $x\in G$ has a resolution $\{Q_{x}^{j}(G):1\leq j \leq \frac {g+s}{2}\}$ into $\mathrm{PPC}$s of $(G\cup S)\setminus \{x\}$.
Each $\mathcal{B}_{x}(G)$ with $x\in S$ has a resolution $\{Q_{x}^{j}(G):1\leq j \leq \frac{g+s}{2}\}$ into $\mathrm{PPC}$s of $(G\cup S)\setminus\{x\}$; particularly, each $Q_{x}^{j}(G)$ with $\frac{g}{2}+1\leq i \leq \frac{g+s}{2}$ is a $\mathrm{PPC}$ of $G$.
Define
$$\mathcal{F} = \mathcal{A}\bigcup\left(\bigcup_{G\in\mathcal{G}}\mathcal{B}(G)\right).$$
It is straightforward that $(X,\mathcal{F})$ forms an $\mathrm{SQS}(gn+s)$.
For $x\in G$ and $G\in \mathcal{G}$, the derived design at $x$ has the block set
$\mathcal{A}_{x}\bigcup \mathcal{B}_{x}(G)$ and a resolution
$$\bigg\{P_{x}^{j}\cup Q_{x}^{j}(G):1\leq j \leq\frac {g+s}{2}\bigg\}\bigcup\bigg\{P_{x}^{j}:\frac {g+s}{2}+1 \leq j \leq \frac {gn+s}{2}\bigg\}.$$
For each $x\in S$, the derived design at $x$ has the block set
$\mathcal{A}_{x}\bigcup (\bigcup_{G\in \mathcal{G}}\mathcal{B}_{x}(G))$ and a resolution
$$\bigg\{P_{x}^{j}(G)\cup Q_{x}^{j}(G):G\in \mathcal{G},1\leq j \leq\frac {g}{2}\bigg\}\bigcup\bigg\{\bigcup_{G\in \mathcal{G}} Q_{x}^{j}(G):\frac {g}{2}+1 \leq j \leq \frac {g+s}{2}\bigg\}.$$
Hence, we end up with an $\mathrm{mcDSQS}(gn+s)$.
\qed

{\Construction \label{C3} Let $mn+s\equiv mk+s\equiv 2,4\pmod 6$. Suppose that there exists an $\mathrm{RDS}(3,k+1,n+1)$. If there exist a $\mathrm{gcDCQS}(m^{k}: s)$ and an $\mathrm{RDGDD}(3,4, (k+1)\{m\})$, then there exists a $\mathrm{gcDCQS}(m^n: s)$.}
\proof See \cite[Lemma 5.5]{Zhou-lkts} for the case $mn+s \equiv mk+s\equiv 4\pmod 6$. Next, we prove the case  $mn+s \equiv mk+s\equiv 2\pmod 6$. Suppose $(I_n \cup\{\infty\}, \mathcal{B})$ is the given $\mathrm{RDS}(3, k+1, n+1)$. Let $\mathcal{B}_\infty=\{B \backslash\{\infty\}: \infty\in B, B \in \mathcal{B}\}$ and $\mathcal{T}=\{B: \infty \notin B, B \in \mathcal{B}\}$. Thus $(I_n, \mathcal{B}_{\infty}, \mathcal{T})$ forms a $1$-fan $\mathrm{S}(3,(\{k\}, \{k+1\}), n)$, where $(I_{n},\mathcal{B}_{\infty})$ is an $\mathrm{S}(2,k,n)$. Let $S=\{\alpha_1, \alpha_2,\ldots, \alpha_s\}, Y=(I_n \times I_m) \cup S$, $\mathcal{H}=\{\{x\} \times I_m: x \in I_n\}$ and $\mathcal{H}_{B} = \{\{x\}\times I_{m}:x\in B\}$, where $B\subseteq I_{n}$.

For any $B\in \mathcal{B}_{\infty}$, we can construct by assumption a $\mathrm{gcDCQS}(m^{k}: s)$ $((B\times I_{m})\cup S,S,\mathcal{H}_{B},\mathcal{A}_B)$; its derived designs are $\mathcal{A}_{B}(x,i)$ with $(x,i)\in B\times I_{m}$ and $\mathcal{A}_{B}(\alpha)$ with $\alpha\in S$.
Moreover, the following hold:\vspace{-0.2cm}
\begin{itemize}
  \item [$\bullet$] each $\mathcal{A}_{B}(x, i)$, where $(x, i) \in B \times I_m$ and $B \in \mathcal{B}_{\infty}$, is the block set of a $\mathrm{gcSTS}(mk+s-1, m+s-1)$ on $((B \times I_m) \setminus\{(x, i)\}) \cup S$ with a hole $(\{x\} \times(I_m \setminus\{i\})) \cup S$;
       its resolution is $\{\mathcal{A}_{B}^{j}(x, i):1\leq j \leq \frac{mk+s}{2}\}$, where each $\mathcal{A}_{B}^{j}(x, i)$ is a $\mathrm{PPC}$ of $((B \times I_m) \setminus\{(x, i)\}) \cup S$; particularly, those $\mathcal{A}_{B}^{j}(x, i)$ with $\frac{m(k-1)}{2} +1 \leq j\leq \frac{mk+s}{2}$ are $\mathrm{PPC}$s of  $(B \setminus\{x\})\times I_{m}$;\vspace{-0.2cm}
  \item [$\bullet$] each $(B\times I_{m},\mathcal{H}_{B},\mathcal{A}_B(\alpha))$, where $\alpha\in S$, is a $\mathrm{KF}(m^k)$; its resolution is $\{\mathcal{A}_{B}^{y,j}(\alpha):y\in B, 1\leq j\leq \frac{m}{2}\}$, where each $\mathcal{A}_{B}^{y,j}(\alpha)$ is a $\mathrm{PC}$ of $(B \setminus\{y\})\times I_{m}$.\vspace{-0.2cm}
\end{itemize}

For any $B\in \mathcal{T}$, we can construct by assumption an $\mathrm{RDGDD}(3,4,(k+1)\{m\})$ $(B\times I_{m},\mathcal{H}_{B},\mathcal{D}_B)$; its derived designs are $\mathcal{D}_{B}(x,i)$ with $(x,i)\in B\times I_{m}$.
Moreover,
each $\mathcal{D}_B(x, i)$ forms the block set of an $\mathrm{RGDD}(2,3, mk)$ with group set $\mathcal{H}_{B\setminus\{x\}}$; its resolution is $\{\mathcal{D}_B^j(x, i): 1 \leqslant j \leqslant \frac{m(k-1)}{2}\}$, each $\mathcal{D}_B^j(x, i)$ being a $\mathrm{PC}$ of $(B \setminus\{x\}) \times I_m$.
Define
$$\mathcal{F} = \left(\bigcup_{B\in \mathcal{B}_{\infty}}\mathcal{A}_{B}\right)\bigcup\left(\bigcup_{B\in \mathcal{T}}\mathcal{D}_{B}\right).$$
Then it is routine to check that $(Y,S,\mathcal{H},\mathcal{F})$ forms a $\mathrm{CQS}(m^{n}:s)$. Next we verify that $(Y,S,\mathcal{H},\mathcal{F})$ is a $\mathrm{gcDCQS}(m^n: s)$.
For $(x,i)\in I_{n}\times I_{m}$, the derived design at $(x,i)$ has the block set
$$\mathcal{F}(x,i) = \left(\bigcup_{x\in B,B\in \mathcal{B}_{\infty}}\mathcal{A}_{B}(x,i)\right)\bigcup\left(\bigcup_{x\in B,B\in \mathcal{T}}\mathcal{D}_{B}(x,i)\right).$$
For $\alpha\in S$, the derived design at $\alpha$ has the block set
$$\mathcal{F}(\alpha)=\bigcup_{B\in \mathcal{B}_{\infty}}\mathcal{A}_{B}(\alpha).$$

For each $x \in I_n\cup \{\infty\}$, let $\mathcal{B}_x=\{B \backslash\{x\}: x \in B, B \in \mathcal{B}\}$. Since $(I_{n}\cup\{\infty\},\mathcal{B})$ is an $\mathrm{RDS}(3,k+1,n+1)$, each $((I_n \cup\{\infty\}) \setminus\{x\}, \mathcal{B}_x)$ forms an $\mathrm{RS}(2, k, n)$. Let $\Gamma_x=\{\mathcal{B}_x(h): 1 \leqslant h \leqslant p\}$, where $p=(n-1) /(k-1)$, be the resolution of $\mathcal{B}_{x}$.
For any $x \in I_n$, $1 \leqslant h \leqslant p$, there is a unique block $M_{x}(h)\in \mathcal{B}_{x}(h)$ containing $\infty$. Therefore, for any $B \in \mathcal{B}_x(h) \setminus\{M_{x}(h)\}$, we must have $B \cup\{x\} \in \mathcal{T}$. Denote $M_{\infty}(h)=(M_{x}(h)\cup\{x\}) \backslash\{\infty\}$. Then $\{M_{\infty}(h) \setminus\{x\}: 1 \leqslant h \leqslant p\}=\{M_{x}(h) \backslash\{\infty\}: 1 \leqslant h \leqslant p\}$ forms a partition of $I_n \backslash\{x\}$.

(1) Consider the coloring of the derived design $\mathcal{F}(x, i),(x, i) \in I_n \times I_m$, which forms the block set of an $\mathrm{STS}(mn+s-1, m+s-1)$ on $Y \setminus\{(x, i)\}$ with a hole $((\{x\} \times I_m)\cup S) \setminus\{(x,i)\}$.

For $1\leq j \leq \frac{m(k-1)}{2}$ and $1\leq h\leq p$, let
$$\mathcal{F}_{h}^{j}(x,i) = \mathcal{A}_{M_{\infty}(h)}^{j}(x,i) \bigcup \left( \bigcup _{B\in \mathcal{B}_{x}(h)\setminus \{M_{x}(h)\}}\mathcal{D}_{B\cup\{x\}}^{j}(x,i)\right).$$

For $\frac{m(k-1)}{2}+1\leq j \leq  \frac{mk+s}{2}$, let
$$\mathcal{F}^{j}(x,i) = \bigcup_{h=1}^{p}\mathcal{A}_{M_{\infty}(h)}^{j}(x,i).$$
It is immediate that the collection $\{\mathcal{F}_{h}^{j}(x,i):1\leq j \leq \frac{m(k-1)}{2},1\leq h\leq p\}\bigcup \{\mathcal{F}^{j}(x,i):\frac{m(k-1)}{2}+1\leq j \leq  \frac{mk+s}{2}\}$ forms a resolution of $\mathcal{F}(x,i)$ into $\mathrm{PPC}$s of $Y\setminus \{(x,i)\}$; in particular, each $\mathcal{F}^{j}(x,i)$ with $\frac{m(k-1)}{2}+1\leq j \leq  \frac{mk+s}{2}$ is a $\mathrm{PPC}$ of $(I_{n}\setminus \{x\} )\times I_m$. As a result, $(Y\setminus\{(x,i)\},((\{x\} \times I_m)\cup S) \setminus\{(x,i)\},\mathcal{F}(x, i))$ is good colorable.

(2) Consider the coloring of the derived design $\mathcal{F}(\alpha)$, $\alpha\in S$, which forms the block set of a $\mathrm{GDD}(2,3,mn)$ of type $m^{n}$ on $I_{n}\times I_{m}$ with group set $\mathcal{H}$.

For $y\in I_{n}$ and $1\leq j\leq \frac{m}{2}$, let
$$\mathcal{F}^{y,j}(\alpha) =  \bigcup_{y\in B,B\in \mathcal{B}_{\infty}}\mathcal{A}_{B}^{y,j}(\alpha).$$
We can check that $\mathcal{F}(\alpha)$ is the block set of a $\mathrm{KF}(m^{n})$ with a resolution $\{\mathcal{F}^{y,j}(\alpha):y\in I_{n},1\leq j\leq \frac{m}{2}\}$,
where each $\mathcal{F}^{y,j}(\alpha)$ is a $\mathrm{PC}$ of $
(I_{n}\setminus\{y\})\times I_{m}$.
Hence, we end up with a $\mathrm{gcDCQS}(m^n: s)$.
\qed

{\Construction \label{C2} Suppose that there exists a $1$-$\mathrm{FG}(3,(K_{1},K),gn)$ of type $g^{n}$.
Let $m,s$ be integers such that $mgn+s\equiv mk_{1}+s\equiv 2,4 \pmod 6$ for any $k_{1}\in K_{1}$.
 If there exist a $\mathrm{gcDCQS}(m^{k_{1}}: s)$ for any $k_{1}\in K_{1}$ and an $\mathrm{FDGDD}(3,4, k\{m\})$ for any $k\in K$, then there exists a $\mathrm{gcDCQS}((mg)^n: s)$.}
\proof
See \cite[Lemma 5.14]{Zhou-lkts} for the case $mgn+s\equiv mk_{1}+s\equiv 4\pmod 6$. Next, we only prove the case  $mgn+s\equiv mk_{1}+s\equiv 2\pmod 6$. 
Suppose $(X,\mathcal{G},\mathcal{B},\mathcal{T})$ is the given $1$-$\mathrm{FG}(3,(K_{1},K),gn)$ of type $g^{n}$. Let $S=\{\alpha_{1},\alpha_{2},\ldots, \alpha_{s}\}$, $Y=(X \times I_m) \cup S$, $\mathcal{H} = \{G\times I_{m}:G\in \mathcal{G}\}$ and $\mathcal{H}_{B} = \{\{x\}\times I_{m}:x\in B\}$, where $B\subseteq X$.

For any $B\in \mathcal{B}$, construct a $\mathrm{gcDCQS}(m^{|B|}: s)$ $((B\times I_{m})\cup S,S, \mathcal{H}_{B}, \mathcal{A}_{B})$; its derived designs are
$\mathcal{A}_{B}(x,i)$ with $(x,i)\in B\times I_{m}$ and $\mathcal{A}_{B}(\alpha)$ with $\alpha\in S$. Moreover, we have the following:\vspace{-0.2cm}
\begin{itemize}
  \item [$\bullet$] each $\mathcal{A}_{B}(x, i)$, where $(x, i) \in B \times I_m$, is the block set of a $\mathrm{gcSTS}(m|B|+s-1, m+s-1)$ on $((B \times I_m) \setminus\{(x, i)\}) \cup S$ with a hole $(\{x\} \times(I_m \setminus\{i\})) \cup S$;
       its resolution is $\{\mathcal{A}_{B}^{y,j}(x, i):y\neq x, y\in B,1\leq j \leq \frac{m}{2}\}\cup\{\mathcal{A}_{B}^{j}(x, i):1 \leq j\leq \frac{m+s}{2}\}$, where each $\mathcal{A}_{B}^{y,j}(x, i)$ is a $\mathrm{PPC}$ of $((B \times I_m) \setminus\{(x, i)\}) \cup S$, and those $\mathcal{A}_{B}^{j}(x, i)$ are $\mathrm{PPC}$s of  $(B \setminus\{x\})\times I_{m}$;\vspace{-0.2cm}

  \item [$\bullet$] each $(B\times I_{m},\mathcal{H}_{B},\mathcal{A}_B(\alpha))$, where $\alpha\in S$, is a $\mathrm{KF}(m^{|B|})$; its resolution is $\{\mathcal{A}_{B}^{y,j}(\alpha):y\in B, 1\leq j\leq \frac{m}{2}\}$, where each $\mathcal{A}_{B}^{y,j}(\alpha)$ is a $\mathrm{PC}$ of $(B \setminus\{y\})\times I_{m}$.\vspace{-0.2cm}
\end{itemize}

For any $B\in \mathcal{T}$, we can construct by assumption an $\mathrm{FDGDD}(3,4,|B|\{m\})$ $(B\times I_{m},\mathcal{H}_{B},\mathcal{D}_B)$; its derived designs are $\mathcal{D}_{B}(x,i)$ with $(x,i)\in B\times I_{m}$.
Moreover, each $((B\setminus\{x\})\times I_{m},\mathcal{H}_{B\setminus\{x\}},$ $\mathcal{D}_{B}(x,i))$ is a $\mathrm{KF}(m^{|B|-1})$; its resolution is $\{\mathcal{D}_{B}^{y,j}(x,i):y\in B\setminus\{x\}, 1\leq j\leq \frac{m}{2}\}$, where each $\mathcal{D}_{B}^{y,j}(x,i)$ is a $\mathrm{PC}$ of $(B \setminus\{x,y\})\times I_{m}$.
Define

$$\mathcal{F} = \left(\bigcup_{B\in \mathcal{B}}\mathcal{A}_{B}\right)\bigcup\left(\bigcup_{B\in \mathcal{T}}\mathcal{D}_{B}\right).$$
Then it is routine to check that $(Y,S,\mathcal{H},\mathcal{F})$ forms a $\mathrm{CQS}((mg)^{n}:s)$. Next we verify that $(Y,S,\mathcal{H},\mathcal{F})$ is a $\mathrm{gcDCQS}((mg)^n: s)$.
For $(x,i)\in X\times I_{m}$, the derived design at $(x,i)$ has the block set
$$\mathcal{F}(x,i) = \left(\bigcup_{x\in B,B\in \mathcal{B}}\mathcal{A}_{B}(x,i)\right)\bigcup\left(\bigcup_{x\in B,B\in \mathcal{T}}\mathcal{D}_{B}(x,i)\right).$$
For $\alpha\in S$, the derived design at $\alpha$ has the block set
$$\mathcal{F}(\alpha)=\bigcup_{B\in \mathcal{B}}\mathcal{A}_{B}(\alpha).$$

(1) For any $x\in G$, $G\in \mathcal{G}$ and $i \in I_m$, consider the coloring of the derived design $\mathcal{F}(x, i)$, which forms the block set of an $\mathrm{STS}(mgn+s-1, mg+s-1)$ on $Y \setminus\{(x, i)\}$ with a hole $((G \times I_m)\cup S) \setminus\{(x,i)\} $.

For any $y\in X\setminus G$ and $1\leq j \leq \frac{m}{2}$, let
$$\mathcal{F}^{y,j}(x,i) = \mathcal{A}_{B_{1}}^{y,j}(x,i) \bigcup \left( \bigcup _{\{x,y\}\subset B,B\in \mathcal{T}}\mathcal{D}_{B}^{y,j}(x,i)\right),$$
where $B_{1}$ is the unique block of $\mathcal{B}$ such that $\{x,y\}\subset B_{1}$.
For any $y\in G\setminus \{x\}$ and $1\leq j \leq \frac{m}{2}$, let
$$\mathcal{F}^{y,j}(x,i) =  \bigcup _{\{x,y\}\subset B,B\in \mathcal{T}}\mathcal{D}_{B}^{y,j}(x,i).$$
For $1\leq j \leq  \frac{m+s}{2}$, let
$$\mathcal{F}^{j}(x,i) = \bigcup_{x\in B, B \in \mathcal{B}}\mathcal{A}_{B}^{j}(x,i).$$
Then it follows that the collection $\{\mathcal{F}^{y,j}(x,i):y\in X\setminus \{x\},1\leq j \leq \frac{m}{2}\}\cup \{\mathcal{F}^{j}(x,i):1\leq j \leq  \frac{s+m}{2}\}$ forms a resolution of $\mathcal{F}(x,i)$, where each $\mathcal{F}^{y,j}(x,i)$ with $y\in X\setminus G,1\leq j \leq \frac{m}{2}$ is a $\mathrm{PPC}$ of $Y\setminus \{(x,i)\}$ and those $\mathcal{F}^{y,j}(x,i)$ with $y\in G\setminus \{x\},1\leq j \leq \frac{m}{2}$ and $\mathcal{F}^{j}(x,i)$ with $1\leq j \leq  \frac{m+s}{2}$ are $\mathrm{PPC}$s of $Y\setminus((G \times I_m) \cup S)$.

(2) Consider the coloring of the derived design $\mathcal{F}(\alpha)$, $\alpha\in S$, which forms the block set of a $\mathrm{GDD}(2,3,mgn)$ of type $(mg)^{n}$ on $X\times I_{m}$ with group set $\mathcal{H}$.
It is obvious that it is a $\mathrm{KF}((mg)^{n})$ with a resolution $\{\mathcal{F}^{y,j}(\alpha):y\in X,1\leq j\leq \frac{m}{2}\}$, where
$$\mathcal{F}^{y,j}(\alpha) =  \bigcup_{y\in B,B\in \mathcal{B}}\mathcal{A}_{B}^{y,j}(\alpha).$$
Hence, we end up with a $\mathrm{gcDCQS}((mg)^n: s)$.
\qed

{\noindent\textbf{Remark:}} From the proof of Construction \ref{C2}, for the case of $mgn+s\equiv 2\pmod 6$, one can weaken the restrictions to $mgn+s\equiv mk_{1}+s\equiv 2\pmod 6$ for some $k_{1}\in K_{1}$.

The rest of this subsection introduces a special combinatorial structure denoted $2$-$\mathrm{RSQS}^{*}(v)$, which is proposed in \cite{DRSQS11}. This design will be used to present another recursive construction for $\mathrm{gcDCQS}$s.

Let $X$ be a set of $v$ elements. A collection $\{(X,\mathcal{B}_{k}^{l}):1\leq k \leq \frac{v-1}{3},1\leq l \leq 3\}$ of $v-1$ $\mathrm{S}(2,4,v)$s is denoted by $2$-$\mathrm{RSQS}^{*}(v)$ if this collection has the following properties:\vspace{-0.2cm}
\begin{itemize}
  \item [(i)] for $1\leq k \leq \frac{v-1}{3}$, the three block sets $\mathcal{B}_{k}^{l}$ $(1\leq l \leq 3)$ share a common $\mathrm{PC}$ of $X$, say $P_{k}$;\vspace{-0.2cm}
  \item [(ii)] $(X,\mathcal{B}^{'})$, where $\mathcal{B}^{'} = \bigcup_{1\leq k \leq \frac{v-1}{3}}P_{k}$, is an $\mathrm{S}(2,4,v)$ (called the special $\mathrm{RS}(2,4,v)$);\vspace{-0.2cm}
  \item [(iii)] $(X,\mathcal{B})$, where $\mathcal{B}$ is the usual union of $\mathcal{B}_{k}^{l}$ ($1\leq k \leq \frac{v-1}{3}$, $1\leq l \leq 3$) is an $\mathrm{SQS}(v)$ (called the underlying $\mathrm{SQS}(v)$);\vspace{-0.2cm}
  \item [(iv)] each quadruple of $\mathcal{B}\setminus \mathcal{B}^{'}$ occurs twice in the multiset $\bigcup_{1\leq k \leq \frac{v-1}{3},1\leq l \leq 3}\mathcal{B}_{k}^{l}$, while each quadruple of $\mathcal{B}^{'}$ occurs three times.
\end{itemize}

{\Construction \label{C4} Let $4m+s\equiv 2,4\pmod 6$. Suppose that there exists a $2$-$\mathrm{RSQS}^{*}(v)$. If there exist a $\mathrm{gcDCQS}(m^{4}:s)$ and an $\mathrm{RDTD}(3,4,m)$, then there exists a $\mathrm{gcDCQS}(m^{v}:s)$.}
\proof Note that $v\equiv 4\pmod {12}$ is the necessary condition for the existence of a $2$-$\mathrm{RSQS}^{*}(v)$. Then one readily checks that $mv+s\equiv4m+s\equiv 2,4\pmod 6$. Let the collection $\{(X,\mathcal{B}_{k}^{l}):1\leq k \leq \frac{v-1}{3},1\leq l \leq 3\}$ be the given $2$-$\mathrm{RSQS}^{*}(v)$.
Let $(X,\mathcal{B})$ be the underlying $\mathrm{SQS}(v)$ and $\mathcal{B}^{'}$ be the spacial $\mathrm{RS}(2,4,v)$.
Now, we consider the derived designs of this $2$-$\mathrm{RSQS}^{*}(v)$. For any $x\in X$,  we denote the derived design at $x$ of $\mathcal{B}$, $\mathcal{B}_{k}^{l}$ and $P_{k}$ by $\mathcal{B}_{x}$, $\mathcal{B}_{x,k}^{l}$ and $P_{x,k}$, respectively. They have the following properties:\vspace{-0.2cm}
\begin{itemize}
  \item [$\bullet$] for $1\leq k \leq \frac{v-1}{3}$, the three block sets $\mathcal{B}_{x,k}^{l}$ $(1\leq l \leq 3)$ share a common block $P_{x,k}$;\vspace{-0.2cm}
  \item [$\bullet$] $\mathcal{B}_{x}^{'}$, where $\mathcal{B}_{x}^{'} = \bigcup_{1\leq k \leq \frac{v-1}{3}}P_{x,k}$, is a $\mathrm{PC}$ of $X\setminus\{x\}$;\vspace{-0.2cm}
  \item [$\bullet$] each triple of $\mathcal{B}_{x}\setminus \mathcal{B}_{x}^{'}$ occurs twice in the multiset $\bigcup_{1\leq k \leq \frac{v-1}{3},1\leq l \leq 3}\mathcal{B}_{x,k}^{l}$, while each triple of $\mathcal{B}_{x}^{'}$ occurs three times.\vspace{-0.2cm}
\end{itemize}

Let $S=\{\alpha_{1},\alpha_{2},\ldots,\alpha_{s}\}$, $Y = (X\times I_{m})\cup S$ and $\mathcal{G}_{B} =\{\{x\}\times I_{m}:x\in B\}$, where $B\subseteq X$. We will construct the desired $\mathrm{gcDCQS}(m^{v}:s)$ on the point set $Y$ with stem $S$ and group set $\mathcal{G}_{X}$.

For each block $B\in \mathcal{B}\setminus\mathcal{B}^{'}$, we construct an $\mathrm{RDTD}(3,4,m)$ $(B\times I_{m}, \mathcal{G}_{B},\mathcal{A}_{B})$. For every point $(x,i)\in B\times I_{m}$, denote the corresponding derived design as $((B\setminus\{x\})\times I_{m}, \mathcal{G}_{B\setminus\{x\}},\mathcal{A}_{B}(x,i))$ and each $\mathcal{A}_{B}(x,i)$ has a resolution $\{\mathcal{A}_{B}^{j}(x,i):1\leq j \leq m\}$, where each $\mathcal{A}_{B}^{j}(x,i)$ is a $\mathrm{PC}$ of $(B\setminus \{x\})\times I_{m}$.

For each block $B\in\mathcal{B}^{'}$, we construct a $\mathrm{gcDCQS}(m^{4}:s)$ $((B\times I_{m})\cup S, S, \mathcal{G}_{B},\mathcal{C}_{B})$. Let
\begin{align*}
\delta = \begin{cases}
\frac{m+s}{2}, & \text{ if } 4m+s\equiv  2 \pmod{6},\\
\frac{m+s}{2}-1, & \text{ if } 4m+s\equiv 4 \pmod{6}.\\
		   \end{cases}
\end{align*}
For every point $(x,i)\in B\times I_{m}$, denote the corresponding derived design as $(((B\times I_{m})\cup S)\setminus\{(x,i)\}, ((\{x\}\times I_{m})\cup S)\setminus\{x,i\},\mathcal{C}_{B}(x,i))$ and each $\mathcal{C}_{B}(x,i)$ has a resolution $\{\mathcal{C}_{B}^{j}(x,i):1\leq j \leq \frac{3m}{2}+\delta\}$, where each $\mathcal{C}_{B}^{j}(x,i)$ is a $\mathrm{PPC}$ of $((B\times I_{m})\cup S)\setminus\{(x,i)\}$; particularly, those $\mathcal{C}_{B}^{j}(x,i)$ with $\frac{3m}{2} + 1 \leq j \leq \frac{3m}{2}+\delta$ are $\mathrm{PPC}$s of $(B\setminus \{x\})\times I_{m}$.
For every point $\alpha\in S$, denote the corresponding derived design as $(B\times I_{m},\mathcal{G}_{B},\mathcal{C}_{B}(\alpha))$ and each $\mathcal{C}_{B}(\alpha)$ has a resolution $\{\mathcal{C}_{B}^{\eta,j}(\alpha):1\leq j \leq \frac{m}{2},\eta\in B\}$, where each $\mathcal{C}_{B}^{\eta,j}(\alpha)$ is a $\mathrm{PC}$ of $(B\setminus\{\eta\})\times I_{m}$.

Let
$$\mathcal{F} = \bigg(\bigcup_{B\in \mathcal{B}\setminus\mathcal{B}^{'}}\mathcal{A}_{B}\bigg)\bigcup \bigg(\bigcup_{B\in \mathcal{B}^{'}}\mathcal{C}_{B}\bigg).$$
Clearly, $(Y,S,\mathcal{G}_{X},\mathcal{F})$ is the block set of a $\mathrm{CQS}(m^{v}:s)$. For $(x,i)\in X\times I_{m}$, the derived design at $(x,i)$ has the block set
$$\mathcal{F}_{(x,i)} = \bigg(\bigcup_{x\in B,B\in \mathcal{B}\setminus\mathcal{B}^{'}}\mathcal{A}_{B}(x,i)\bigg)\bigcup \bigg(\bigcup_{x\in B,B\in \mathcal{B}^{'}}\mathcal{C}_{B}(x,i)\bigg).$$
Next, we check that each block set $\mathcal{F}_{(x,i)}$ can be partitioned into $\frac{m(v-1)}{2}$ $\mathrm{PPC}$s of $Y\setminus\{(x,i)\}$ and $\delta$ $\mathrm{PPC}$s of $(X\setminus\{x\})\times I_{m}$.

For $1\leq k \leq \frac{v-1}{3}$, $1\leq l \leq 3$ and $1\leq r \leq \frac{m}{2}$, let
$$\mathcal{F}_{(x,i)}(k,l,r) = \bigg(\bigcup_{x\in B,B\setminus\{x\}\in \mathcal{B}_{x,k}^{l}\setminus \{P_{x,k}\}}\mathcal{A}_{B}^{r^{\prime}}(x,i)\bigg)\bigcup \mathcal{C}_{P_{x,k}\cup\{x\}}^{\frac{m}{2}(l-1)+r}(x,i),$$
where the value $r^{\prime}$ equals $r$ or $r+\frac{m}{2}$ according to the first or the second occurrence of $B$ in $\mathcal{B}\setminus\mathcal{B}^{'}$.

For $1 \leq j \leq \delta$, let
$$\mathcal{F}_{(x,i)}(j) = \bigcup_{P_{x,k}\in \mathcal{B}_{x}^{'},1\leq k \leq \frac{v-1}{3}}\mathcal{C}_{P_{x,k}\cup\{x\}}^{\frac{3m}{2}+j}(x,i).$$
So each $\mathcal{F}_{(x,i)}$ has a  desired resolution $\{\mathcal{F}_{(x,i)}(k,l,r):1\leq k \leq \frac{v-1}{3},1\leq l \leq 3,1\leq r \leq \frac{m}{2}\}\bigcup\{\mathcal{F}_{(x,i)}(j):1 \leq j \leq \delta\}$.

For $\alpha\in S$, the derived design at $\alpha$ has the block set
$$\mathcal{F}_{\alpha}=\bigcup_{B\in \mathcal{B}^{\prime}}\mathcal{C}_{B}(\alpha).$$
For $\eta\in X$, $1 \leq j \leq \frac{m}{2}$, let
$$\mathcal{F}_{\alpha}(\eta,j)=\bigcup_{\eta\in B,B\in \mathcal{B}^{\prime}}\mathcal{C}_{B}^{\eta,j}(\alpha).$$
So each $\mathcal{F}_{\alpha}$ is a block set of a $\mathrm{KF}(m^{v})$.
Then we obtain a $\mathrm{gcDCQS}(m^{v}:s)$.
\qed

\begin{Theorem}\label{2RDSQS}
\begin{itemize}
  \item [$(1)$] {\cite[Corollary 3.2]{DRSQS11}} There is a $2$-$\mathrm{RSQS}^{*}(4^{n})$ for any integer $n\geq 2$.
 \item [$(2)$] {\cite{J_TD,DRSQS11}} There is an $\mathrm{RDTD}(3,4,m)$ for any positive integer $m$ satisfying $\gcd(m,4) \neq 2$ and $\gcd(m,18) \neq 3$.
 \end{itemize}
\end{Theorem}

\begin{Corollary}\label{Rmul}
Let $m\equiv 4,8\pmod {12}$. If there exist a $\mathrm{gcDCQS}(m^{4}:0)$ and an $\mathrm{mcDSQS}(m)$, then there exists an $\mathrm{mcDSQS}(m\cdot 4^{n})$ for any integer $n\geq 2$.
\end{Corollary}
\proof Let $m\equiv 4,8\pmod {12}$ and $s=0$. For any integer $n\geq 2$, apply Construction \ref{C4}. Since a $2$-$\mathrm{RSQS}^{*}(4^{n})$ and an $\mathrm{RDTD}(3,4,m)$ exist by Theorem \ref{2RDSQS}, if there exists a $\mathrm{gcDCQS}(m^{4}:0)$, then there exists a $\mathrm{gcDCQS}(m^{4^{n}}:0)$. Further, if there exists an $\mathrm{mcDSQS}(m)$,
then there exists an $\mathrm{mcDSQS}(m\cdot 4^{n})$ by Construction \ref{mcdcqs-mcdsqs}.
\qed

\noindent{\bf Remark:} Construction \ref{C4} presents a new construction even for $\mathrm{RDCQS}$s. As a special case, taking $m=4$ and $s=0$ in Construction \ref{C4} gives the quadrupling construction by Lui and Lei \cite[Construction 3.2]{LEI1} (noting the existence of an RDCQS$(4^4:0)$ and an $\mathrm{RDTD}(3,4,4)$).

\subsection{Special considerations}

Construction \ref{mcdcqs-mcdsqs} provides a construction of an $\mathrm{mcDSQS}(gn+s)$ by filling in holes of a $\mathrm{gcDCQS}(g^{n}:s)$. Constructions \ref{C3}-\ref{C4} present approaches of producing $\mathrm{gcDCQS}$s. However, in order to deal with $\mathrm{mcDSQS}(v)$ where $v\equiv 2\pmod 6$, we have to construct some input designs, say $\mathrm{gcDCQS}(g^{k}:s)$ with $g+s\geq 20$ (because $\chi^{\prime}(v) >\frac{v+1}{2}$ if $v\in\{7,13\}$), which requires strong computing power.
In this subsection, we consider $\mathrm{CQS}(g^{n}:s)$ with group size $g \equiv2\pmod 6$ and stem size $s=2$. We introduce two new types of colorings of the derived designs and define $\mathrm{^{1}cDCQS}(g^{n}:2)$ and $\mathrm{^{2}cDCQS}(g^{n}:2)$; $\mathrm{^{1}cDCQS}$ is used in a new filling construction 
and $\mathrm{^{2}cDCQS}$ is an instrumental design to generate $\mathrm{^{1}cDCQS}$s in a recursive construction.

Let $(X,\mathcal{G},\mathcal{B})$ be a $\mathrm{GDD}(2,3,gn)$ of type $g^{n}$ with $g$ even. We say that $(X,\mathcal{G},\mathcal{B})$ is {\em{good colorable}} and denoted by $\mathrm{gcGDD}(2,3,gn)$ of type $g^{n}$ if the block set $\mathcal{B}$ can be partitioned into $\frac{gn}{2}+1$ $\mathrm{PPC}$s of $X\setminus G$, $G\in \mathcal{G}$, to be specific,  $\mathcal{B} = \left(\bigcup_{G\in \mathcal{G},1\leq i \leq \frac{g}{2}} P_{G}^{i}\right)\bigcup P_{H}^{*}$, where $P_{G}^{i}$ is a $\mathrm{PPC}$ of $X \setminus G$ and $P_{H}^{*}$ is a $\mathrm{PPC}$ of $X\setminus H$ for some $H\in \mathcal{G}$. We call $H$ the special group of the $\mathrm{gcGDD}(2,3,gn)$ of type $g^{n}$.

Now, we define two coloring types for the derived designs of a $\mathrm{CQS}(g^{n}:2)$. Let $(X, S, \mathcal{G}, \mathcal{A})$ be a $\mathrm{CQS}(g^n: 2)$. For $x \in X$, put $\mathcal{A}_x=\{A \setminus\{x\}: x \in A, A \in \mathcal{A}\}$.
A $\mathrm{CQS}(g^n: 2)$ whose derived designs admit type I colorings, denoted by $^{1}\mathrm{cDCQS}(g^n: 2)$, refers to a $\mathrm{CQS}(g^n: 2)$ whose derived designs $\mathcal{A}_x$, $x \in X$, have the following properties:\vspace{-0.2cm}
\begin{itemize}
  \item [(1)] for $x\in G$ and $G\in \mathcal{G}$, each $(X\setminus \{x\}, (G\cup S)\setminus \{x\},\mathcal{A}_{x})$ is a $\mathrm{gcSTS}(gn+1,g+1)$;\vspace{-0.2cm}
  \item [(2)] for $x\in S$, each $(X\setminus S, \mathcal{G}, \mathcal{A}_{x})$ is a $\mathrm{gcGDD}(2,3,gn)$ of type $g^{n}$ with the same special group.\vspace{-0.2cm}
\end{itemize}
A $\mathrm{CQS}(g^n: 2)$ whose derived designs admit type II colorings, denoted by $^{2}\mathrm{cDCQS}(g^n: 2)$, refers to a $\mathrm{CQS}(g^n: 2)$ whose derived designs $\mathcal{A}_x$, $x \in X$, have the following properties:\vspace{-0.2cm}
\begin{itemize}
  \item [(1)] for $x\in G$ and $G\in \mathcal{G}$, each $(X\setminus \{x\}, (G\cup S)\setminus \{x\},\mathcal{A}_{x})$ is a $\mathrm{gcSTS}(gn+1,g+1)$;\vspace{-0.2cm}

  \item [(2)] for $x\in S$, each $(X\setminus S, \mathcal{G}, \mathcal{A}_{x})$ is an $\mathrm{RGDD}(2,3,gn)$ of type $g^{n}$.
\end{itemize}

{\Construction \label{cqs-wrdsqs1} Let $n\equiv 0\pmod 3$. If there exist a $\mathrm{^{1}cDCQS}(g^n: 2)$ and an $\mathrm{RDSQS}(g+2)$, then there exists an $\mathrm{mcDSQS}(gn+2)$.}
\proof Suppose $(X,S,\mathcal{G},\mathcal{A})$ is the given $^{1}\mathrm{cDCQS}(g^n: 2)$.
For $x \in X$, let $\mathcal{A}_x=\{A \setminus\{x\}: x \in A, A \in \mathcal{A}\}$.
Then, each $(X\setminus \{x\}, (G\cup S)\setminus \{x\},\mathcal{A}_{x})$ with $x\in G$ and $G\in \mathcal{G}$ is a $\mathrm{gcSTS}(gn+1,g+1)$ with a hole $(G\cup S)\setminus\{x\}$; its resolution is $\{P_{x}^{j}:1 \leq j \leq \frac{gn}{2}+1\}$, where each $P_{x}^{j}$ is a $\mathrm{PPC}$ of $X\setminus \{x\}$; particularly, those $P_{x}^{j}$ with $1 \leq j \leq \frac{g}{2}+1$ are $\mathrm{PPC}$s of $X\setminus (G\cup S)$.
In addition, each $(X\setminus S, \mathcal{G}, \mathcal{A}_{x})$ with $x\in S$ is a $\mathrm{gcGDD}(2,3,gn)$ of type $g^{n}$ with the special group $H\in \mathcal{G}$; its resolution is $\{P_{x}^{j}(G):G\in \mathcal{G},1\leq j \leq \frac{g}{2}\}\cup\{P_{x}^{*}\}$, where each $P_{x}^{j}(G)$, $G\in \mathcal{G}$ and $1\leq j \leq \frac{g}{2}$, is a $\mathrm{PPC}$ of $X\setminus (G\cup S)$ and $P_{x}^{*}$ is a $\mathrm{PPC}$ of $X\setminus (H\cup S)$.

For each $G\in \mathcal{G}$, we can construct by assumption an $\mathrm{RDSQS}(g+2)$ on $G\cup S$ with block set $\mathcal{B}(G)$.
So we have $g+2$ derived $\mathrm{KTS}(g+1)$s, say $((G\cup S) \setminus\{x\}, \mathcal{B}_{x}(G))$ with $x\in G\cup S$.
Further, each $\mathcal{B}_{x}(G)$ with $x\in G\cup S$ has a resolution $\{Q_{x}^{j}(G):1\leq j \leq \frac {g}{2}\}$ into $\mathrm{PC}$s of $(G\cup S)\setminus\{x\}$.
Define
$$\mathcal{F} = \mathcal{A}\bigcup\left(\bigcup_{G\in\mathcal{G}}\mathcal{B}(G)\right).$$
Then, $\mathcal{F}$ is the block set of an $\mathrm{SQS}(gn+2)$.
For $x\in G$ and $G\in \mathcal{G}$, the derived design at $x$ has the block set $\mathcal{A}_{x}\cup \mathcal{B}_{x}(G)$ and a resolution
$$\big\{P_{x}^{j}\cup Q_{x}^{j}(G):1\leq j \leq \frac{g}{2}\big\}\bigcup\big\{P_{x}^{j}:\frac{g}{2}+1 \leq j \leq  \frac{gn}{2}+1\big\}.$$
For $x\in S$, the derived design at $x$ has the block set $\mathcal{A}_{x}\bigcup (\bigcup_{G\in \mathcal{G}}\mathcal{B}_{x}(G))$ and a resolution
$$\big\{P_{x}^{j}(G)\cup Q_{x}^{j}(G):G\in \mathcal{G},1\leq j \leq \frac{g}{2}\big\}\bigcup\big\{P_{x}^{*}\big\}.$$
Hence, we end up with the desired $\mathrm{mcDSQS}(gn+2)$.
\qed

{\Construction \label{C1} Let $m\equiv 2 \pmod 6$. Suppose that there exists an $\mathrm{RDS}(3,k+1,n+1)$. If there exist a $\mathrm{^{1}cDCQS}(m^{k}: 2)$, a $\mathrm{^{2}cDCQS}(m^{k}: 2)$ and an $\mathrm{RDGDD}(3,4, (k+1)\{m\})$, then there exists a $\mathrm{^{1}cDCQS}(m^n: 2)$.}

\proof Suppose $(I_n \cup\{\infty\}, \mathcal{B})$ is the given $\mathrm{RDS}(3, k+1, n+1)$. Let $\mathcal{B}_\infty=\{B \backslash\{\infty\}: \infty\in B, B \in \mathcal{B}\}$ and $\mathcal{T}=\{B: \infty \notin B, B \in \mathcal{B}\}$. Thus $(I_n, \mathcal{B}_{\infty},\mathcal{T})$ forms a $1$-fan $\mathrm{S}(3,(\{k\}, \{k+1\}), n)$, where $(I_{n},\mathcal{B}_{\infty})$ is an $\mathrm{S}(2,k,n)$. Let $S=\{\alpha_1, \alpha_2\}, Y=(I_n \times I_m) \cup S$, $\mathcal{H}=\{\{x\} \times I_m: x \in I_n\}$ and $\mathcal{H}_{B} = \{\{x\}\times I_{m}:x\in B\}$, where $B\subseteq I_{n}$.

Fix an element $\eta\in I_{n}$. For any $B\in \mathcal{B}_{\infty}$ and $\eta \in B$, we can construct by assumption a $^{1}\mathrm{cDCQS}(m^{k}: 2)$ $((B\times I_{m})\cup S,S,\mathcal{H}_{B},\mathcal{A}_B)$ with derived designs $\mathcal{A}_{B}(x,i)$ with $(x,i)\in B\times I_{m}$ and $\mathcal{A}_{B}(\alpha)$ with $\alpha\in S$. Moreover, the following should be satisfied:\vspace{-0.2cm}
\begin{itemize}
  \item [(A1)] each $\mathcal{A}_{B}(x, i)$, where $(x, i) \in B \times I_m$, is the block set of a $\mathrm{gcSTS}(mk+1, m+1)$ on $((B \times I_m) \setminus\{(x, i)\}) \cup S$ with a hole $((\{x\} \times I_m)\cup S) \setminus\{(x,i)\}$;
       its resolution is $\{\mathcal{A}_{B}^{j}(x, i):1\leq j \leq \frac{mk}{2}+1\}$, where each $\mathcal{A}_{B}^{j}(x, i)$ is a $\mathrm{PPC}$ of $((B \times I_m) \setminus\{(x, i)\}) \cup S$. Particularly,
       those $\mathcal{A}_{B}^{j}(x, i)$ with $\frac{m(k-1)}{2} +1 \leq j\leq \frac{mk}{2}+1$ are $\mathrm{PPC}$s of  $(B \setminus\{x\})\times I_{m}$;\vspace{-0.2cm}

  \item [(A2)] each $(B\times I_{m},\mathcal{H}_{B},\mathcal{A}_B(\alpha))$, where $\alpha\in S$, is a $\mathrm{gcGDD}(2,3,mk)$ of type $m^{k}$ with the special group $\{\eta\}\times I_{m}$; and $\mathcal{A}_B(\alpha)$ has a resolution $\{\mathcal{A}_{B}^{y,j}(\alpha):y\in B, 1\leq j\leq \frac{m}{2}\}\cup\{\mathcal{A}_{B}^{\eta,0}(\alpha)\}$, where each $\mathcal{A}_{B}^{y,j}(\alpha)$ is a $\mathrm{PPC}$ of $(B \setminus\{y\})\times I_{m}$;\vspace{-0.2cm}
\end{itemize}

For any $B\in \mathcal{B}_{\infty}$ and $\eta \notin B$, we can construct by assumption a $\mathrm{^{2}cDCQS}(m^{k}: 2)$ $((B\times I_{m})\cup S,S,\mathcal{H}_{B},\mathcal{A}_B)$; its derived designs are $\mathcal{A}_{B}(x,i)$ with $(x,i)\in B\times I_{m}$ and $\mathcal{A}_{B}(\alpha)$ with $\alpha\in S$. Moreover, each $\mathcal{A}_{B}(x, i)$, $(x, i) \in B \times I_m$, satisfies the property (A$1$). In addition, each $(B\times I_{m},\mathcal{H}_{B},\mathcal{A}_B(\alpha))$ with $\alpha\in S$ is an $\mathrm{RGDD}(2,3,mk)$ of type $m^{k}$ with a resolution $\{\mathcal{A}_{B}^{j}(\alpha):1\leq j \leq \frac{m(k-1)}{2}\}$ into $\mathrm{PC}$s of $B\times I_{m}$.

For any $B\in \mathcal{T}$, we can construct by assumption an $\mathrm{RDGDD}(3,4,(k+1)\{m\})$ $(B\times I_{m},\mathcal{H}_{B},\mathcal{D}_B)$; its derived designs are $\mathcal{D}_{B}(x,i)$ with $(x,i)\in B\times I_{m}$.
Moreover, each $\mathcal{D}_B(x, i)$, where $(x, i) \in B \times I_m$ and $B \in \mathcal{T}$, forms the block set of an $\mathrm{RGDD}(2,3, mk)$ with group set $\mathcal{H}_{B\setminus\{x\}}$; its resolution is $\{\mathcal{D}_B^j(x, i): 1 \leqslant j \leqslant \frac{m(k-1)}{2}\}$, each $\mathcal{D}_B^j(x, i)$ being a $\mathrm{PC}$ of $(B \setminus\{x\}) \times I_m$.
Define
$$\mathcal{F} = \left(\bigcup_{B\in \mathcal{B}_{\infty}}\mathcal{A}_{B}\right)\bigcup\left(\bigcup_{B\in \mathcal{T}}\mathcal{D}_{B}\right).$$
Then $(Y,S,\mathcal{H},\mathcal{F})$ forms a $\mathrm{CQS}(m^{n}:2)$. For $(x,i)\in I_{n}\times I_{m}$, the derived design at $(x,i)$ has the block set
$$\mathcal{F}(x,i) = \left(\bigcup_{x\in B,B\in \mathcal{B}_{\infty}}\mathcal{A}_{B}(x,i)\right)\bigcup\left(\bigcup_{x\in B,B\in \mathcal{T}}\mathcal{D}_{B}(x,i)\right).$$
For $\alpha\in S$, the derived design at $\alpha$ has the block set
$$\mathcal{F}(\alpha)=\bigcup_{B\in \mathcal{B}_{\infty}}\mathcal{A}_{B}(\alpha).$$
By the proof of Construction \ref{C3}, each derived design $\mathcal{F}(x, i)$, $(x, i) \in I_n \times I_m$, forms the block set of a $\mathrm{gcSTS}(mn+1, m+1)$ on $Y \setminus\{(x, i)\}$ with a hole $(\{x\} \times(I_m \setminus\{i\})) \cup S$.

Next, consider the coloring of the derived design $\mathcal{F}(\alpha)$, $\alpha\in S$.
Since $\mathcal{B}_{\infty}$ is a resolvable $\mathrm{S}(2,k,n)$, it has a resolution $\{\mathcal{B}_{\infty}(h): 1 \leqslant h \leqslant p\}$, where $p=(n-1) /(k-1)$. So there is a unique block $C_{\infty}(h)$ in $\mathcal{B}_{\infty}(h)$, $1 \leqslant h \leqslant p$, that contains $\eta$.
For any block $C_{\infty}(h)$, $1 \leqslant h \leqslant p$, we define a bijection $f_{B}$ from $(B\setminus\{\eta\})\times \{1,2,\ldots, \frac{m}{2}\}$ to $\{1,2,\ldots, \frac{m(k-1)}{2}\}$ arbitrarily.
The collection $\{C_{\infty}(h)\setminus \{\eta\}:1 \leqslant h \leqslant p\}$ is a partition of $I_{n}\setminus \{\eta\}$.

For every $y\in I_{n}$, $y\neq\eta$ and $1\leq j\leq \frac{m}{2}$, there is a unique $h$ such that $y\in C_{\infty}(h)$. Let
$$\mathcal{F}^{y,j}(\alpha) = \mathcal{A}_{C_{\infty}(h)}^{y,j}(\alpha)\bigcup \left(\bigcup_{B\in \mathcal{B}_{\infty}(h)\setminus \{C_{\infty}(h)\}}\mathcal{A}_{B}^{f_{C_{\infty}(h)}(y,j)}(\alpha)\right).$$
For $0\leq j\leq \frac{m}{2}$, let
$$\mathcal{F}^{\eta,j}(\alpha) = \bigcup_{B\in \mathcal{B}_{\infty},\eta\in B}\mathcal{A}_{B}^{\eta,j}(\alpha).$$
We can check that the collection $\{\mathcal{F}^{y,j}(\alpha):y\in I_{n},1\leq j\leq \frac{m}{2}\}\bigcup\{\mathcal{F}^{\eta,0}(\alpha)\}$ forms a resolution of $\mathcal{F}(\alpha)$.
In particular,
each $\mathcal{F}^{y,j}(\alpha)$ ($y\in I_{n}$, $y\neq \eta$ and $1\leq j\leq \frac{m}{2}$) is a $\mathrm{PPC}$ of $(I_{n}\setminus\{y\})\times I_{m}$ and each $\mathcal{F}^{\eta,j}(\alpha)$ ($0\leq j\leq \frac{m}{2}$) is a $\mathrm{PPC}$ of $(I_{n}\setminus\{\eta\})\times I_{m}$.
Hence, we end up with a $\mathrm{^{1}cDCQS}(m^n: 2)$.
\qed

\section{Existence and applications}

In this section, we produce new infinite families  of  $\mathrm{RDSQS}(6n+4)$s and  $\mathrm{mcDSQS}(6n+2)$s. Then  we show their applications in both
large sets of Kirkman triple systems and diameter perfect constant-weight codes.

{\Lemma \label{CQS} There exists an $\mathrm{RDCQS}(8^{4}:2)$.}
\proof Let $X=\mathbb{Z}_{32}\cup\{a,b\}$ and $\mathcal{G} = \{G_i: 0\leq i \leq 3\}$, where $G_i=\{i, 4+i, 8+i, 12+i, 16+i, 20+i, 24+i, 28+i\}$.
We construct a $\mathrm{CQS}(8^{4}:2)$ on point set $X$ with stem $S=\{a,b\}$ and group set $\mathcal{G}$.
Let $\mathcal{H}$ be the permutation group generated by the  permutation
$(0,1,\ldots,31)(a,b)$.
Our $\mathrm{CQS}(8^{4}:2)$ $(X,S,\mathcal{G},\mathcal{B})$ admits the automorphism group $\mathcal{H}$ and the base blocks are listed as follows:
\begin{center}
\begin{tabular}{l l l l l l l l }
\{0, 1, 2, $a$\}& \{0, 1, 3, 16\}& \{0, 1, 4, 13\}& \{0, 1, 5, 6\}& \{0, 1, 7, 9\}\\
\{0, 1, 8, 23\}& \{0, 1, 10, 25\}& \{0, 1, 11, 12\}& \{0, 1, 14, 19\}& \{0, 1, 15, 18\}\\
\{0, 1, 17, 30\}& \{0, 1, 20, 29\}& \{0, 1, 24, 26\}& \{0, 2, 4, 18\}& \{0, 2, 5, 24\}\\
\{0, 2, 6, 28\}& \{0, 2, 7, 23\}& \{0, 2, 10, 29\}& \{0, 2, 11, 27\}& \{0, 2, 12, 21\}\\
\{0, 2, 13, 22\}& \{0, 2, 14, 20\}& \{0, 2, 17, $b$\}& \{0, 3, 6, $a$\}& \{0, 3, 7, 20\}\\
\{0, 3, 8, 25\}& \{0, 3, 9, 26\}& \{0, 3, 10, 27\}& \{0, 3, 11, 14\}& \{0, 3, 15, 28\}\\
\{0, 4, 9, 27\}& \{0, 4, 11, 15\}& \{0, 4, 14, 22\}& \{0, 5, 10, $a$\}& \{0, 5, 11, 26\}\\
\{0, 5, 12, 17\}& \{0, 6, 12, 22\}& \{0, 6, 13, 24\}& \{0, 6, 14, 25\}& \{0, 6, 19, $b$\}\\
\{0, 7, 14, $a$\}& \{0, 9, 18, $b$\}& \{0, 10, 21, $b$\}\\
\end{tabular}
\end{center}
The resolutions of derived designs at $0$ and $a$ are listed in Appendix A. Thus, we obtain an
$\mathrm{RDCQS}(8^{4}:2)$.
\qed

{\Theorem\label{result_2} There exists an $\mathrm{RDSQS}(2^{2n+1}+2)$ for any nonnegative integer $n$.}
\proof For $n=0$, an $\mathrm{RDSQS}(4)$ exists obviously. For $n=1$, it is well-known that an $\mathrm{SQS}(10)$ is unique up to isomorphism and its derived designs are all resolvable.
For $n=2$, apply Construction \ref{mcdcqs-mcdsqs}. Since an $\mathrm{RDCQS}(8^{4}:2)$ exists by Lemma \ref{CQS},  an $\mathrm{RDSQS}(34)$ exists.

For $n\geq 3$, apply Construction \ref{C4} with a $2$-$\mathrm{RSQS}^{*}(4^{n-1})$, which exists by Theorem \ref{2RDSQS}. Since an $\mathrm{RDCQS}(8^{4}:2)$ and an $\mathrm{RDTD}(3,4,8)$ exist by Lemma \ref{CQS} and Theorem \ref{2RDSQS},  an $\mathrm{RDCQS}(8^{4^{n-1}}:2)$ exists. Further, apply Construction \ref{mcdcqs-mcdsqs} to obtain an $\mathrm{RDSQS}(2^{2n+1}+2)$.
\qed

\begin{Corollary}\label{result_4}
There exists  an $\mathrm{LKTS}(3\cdot 2^{2n+1}+3)$ for  any nonnegative integer $n$.
\end{Corollary}
\proof Apply Theorems \ref{RDS_1}  and \ref{result_2}. 	
\qed

%
%

{\Lemma\label{MIN_Dsqs} There exists an $\mathrm{mcDSQS}(v)$ for $v=20,26,32$.}
\proof In \cite{SHI}, Shi et al. show that the cyclic $\mathrm{SQS}(20)$ constructed in \cite{CSQS} is an $\mathrm{mcDSQS}(20)$.
We construct a cyclic $\mathrm{mcDSQS}(v)$ on $\mathbb{Z}_{v}$ for $v=26,32$; the base blocks and the resolutions are listed in Appendix B.
\qed

{\Lemma\label{Dcqs20_1} There exists a $\mathrm{^{1}cDCQS}(2^{9}:2)$.}
\proof Let $X=I_{18}\cup\{a,b\}$ and $\mathcal{G} = \{\{i,i+9\}: 0\leq i \leq 8\}$.
We construct a $\mathrm{CQS}(2^{9}:2)$ on the point set $X$ with stem $S=\{a,b\}$ and group set $\mathcal{G}$.
Let $\mathcal{H}$ be a permutation group generated by the following permutation:
$$ (0,3,6,9,12,15)(1,4,7,10,13,16)(2,5,8,11,14,17)(a,b).$$
Our $\mathrm{CQS}(2^{9}:2)$ $(X,S,\mathcal{G},\mathcal{B})$ admits the automorphism group $\mathcal{H}$ and the base blocks are listed as follows. (The underlined blocks belong to short orbits of length $3$ while the double-underlined blocks lie in short orbits of length $2$.)
\begin{center}
\begin{tabular}{l l l l l l l l }
\{0, 1, 2, $a$\}& \{0, 1, 3, 6\}& \{0, 1, 4, 15\}& \{0, 1, 5, 13\}& \{0, 1, 7, 9\} \\
\{0, 1, 8, $b$\}& \{0, 1, 10, 14\}& \{0, 1, 11, 16\}& \{0, 1, 12, 17\}& \{0, 2, 3, 12\}\\
\{0, 2, 4, 5\}& \{0, 2, 6, 15\}& \{0, 2, 7, 17\}& \{0, 2, 8, 16\}& \{0, 2, 10, $b$\}\\
\{0, 2, 13, 14\}& \{0, 3, 8, 14\}& \{0, 3, 10, 13\}& \{0, 3, 11, $b$\}&\{0, 3, 17, $a$\}\\
\{0, 4, 6, 10\}& \{0, 4, 8, 11\}& \{0, 4, 14, $a$\}&\{0, 4, 17, $b$\}& \{0, 5, 8, 12\}\\
\{0, 5, 10, $a$\}& \{0, 5, 14, 17\}& \{0, 5, 16, $b$\}&\{0, 6, 13, $b$\}& \{0, 6, 16, 17\}\\
\{0, 7, 8, 13\}& \{0, 7, 11, $a$\}&\{0, 7, 14, 16\}& \{0, 8, 10, 17\}& \{0, 11, 13, 17\}\\
\{0, 13, 16, $a$\}&\{1, 2, 4, $b$\}&\{1, 2, 5, 14\}& \{1, 2, 13, 16\}&\{1, 4, 11, 17\}\\
\{1, 5, 7, $b$\}& \{1, 7, 14, 17\}&\{2, 5, 8, $b$\}&\underline{\{0, 2, 9, 11\}}& \underline{\{0, 4, 9, 13\}}\\
\underline{\{1, 2, 10, 11\}}& \underline{\{1, 4, 10, 13\}}&\underline{\underline{\{0, 6, 12, $a$\}}}& \underline{\underline{\{1, 7, 13, $a$\}}}& \underline{\underline{\{2, 8, 14, $a$\}}}
\end{tabular}
\end{center}
Then we only need to list the resolution of the derived designs at $0,1,2,a$ and $b$ to check that this $\mathrm{CQS}(2^{9}:2)$ is a $\mathrm{^{1}cDCQS}(2^{9}:2)$. We list below the block set of the derived designs at the point $x\in\{0,1,2\}$; the blocks in each row form a $\mathrm{PPC}$ of $X\setminus\{x\}$ and the last two rows form $\mathrm{PPC}$s of $I_{18}\setminus G$, where $x\in G$, $G\in \mathcal{G}$.
{\small\begin{center}
\begin{tabular}{l l l l l l l l }
$\mathcal{B}_{0}$:
    & \{1, 4, 15\}& \{2, 9, 11\}& \{3, 17, $a$\}& \{5, 8, 12\}& \{6, 13, $b$\}& \{7, 14, 16\}\\
    & \{1, 5, 13\}& \{2, 10, $b$\}& \{3, 4, 7\}& \{6, 16, 17\}& \{8, 15, $a$\}& \{9, 12, 14\}\\
    & \{1, 8, $b$\}& \{2, 3, 12\}& \{4, 14, $a$\}& \{5, 6, 7\}& \{9, 10, 16\}& \{11, 13, 17\}\\
    & \{2, 8, 16\}& \{3, 11, $b$\}& \{4, 9, 13\}& \{5, 14, 17\}& \{6, 12, $a$\}& \{7, 10, 15\}\\
    & \{1, 10, 14\}& \{2, 6, 15\}& \{3, 5, 9\}& \{4, 8, 11\}& \{7, 12, $b$\}& \{13, 16, $a$\}\\
    & \{1, 12, 17\}& \{2, 4, 5\}& \{3, 10, 13\}& \{6, 8, 9\}& \{7, 11, $a$\}& \{14, 15, $b$\}\\
    & \{1, 7, 9\}& \{3, 15, 16\}& \{4, 17, $b$\}& \{5, 10, $a$\}& \{6, 11, 14\}\\
    & \{1, 2, $a$\}& \{5, 16, $b$\}& \{7, 8, 13\}& \{9, 15, 17\}& \{10, 11, 12\}\\
    & \{1, 11, 16\}& \{2, 7, 17\}& \{3, 8, 14\}& \{4, 6, 10\}& \{12, 13, 15\}\\
    & \{1, 3, 6\}& \{2, 13, 14\}& \{4, 12, 16\}& \{5, 11, 15\}& \{8, 10, 17\}\\
\end{tabular}
\end{center}

\begin{center}
\begin{tabular}{l l l l l l l l }
$\mathcal{B}_{1}$: & \{0, 2, $a$\}& \{3, 12, 13\}& \{4, 11, 17\}& \{5, 7, $b$\}& \{6, 10, 15\}& \{8, 14, 16\}\\
    & \{0, 3, 6\}& \{2, 7, 12\}& \{4, 10, 13\}& \{5, 8, 15\}& \{9, 14, $b$\}& \{16, 17, $a$\}\\
    & \{0, 4, 15\}& \{2, 10, 11\}& \{3, 16, $b$\}& \{5, 12, $a$\}& \{6, 13, 14\}& \{8, 9, 17\}\\
    & \{0, 8, $b$\}& \{2, 13, 16\}& \{3, 10, 17\}& \{4, 9, 12\}& \{6, 7, 11\}& \{14, 15, $a$\}\\
    & \{2, 4, $b$\}& \{3, 8, $a$\}& \{5, 6, 17\}& \{7, 10, 16\}& \{9, 13, 15\}& \{11, 12, 14\}\\
    & \{0, 10, 14\}& \{2, 6, 8\}& \{3, 7, 15\}& \{4, 5, 16\}& \{9, 11, $a$\}& \{13, 17, $b$\}\\
    & \{0, 11, 16\}& \{2, 15, 17\}& \{3, 4, 14\}& \{5, 9, 10\}& \{6, 12, $b$\}& \{7, 13, $a$\}\\
    & \{0, 5, 13\}& \{2, 3, 9\}& \{4, 6, $a$\}& \{7, 14, 17\}& \{8, 10, 12\}& \{11, 15, $b$\}\\
    & \{0, 7, 9\}& \{2, 5, 14\}& \{8, 11, 13\}&  \{12, 15, 16\}\\
    & \{0, 12, 17\}& \{3, 5, 11\}& \{4, 7, 8\}& \{6, 9, 16\}\\
\end{tabular}
\end{center}

\begin{center}
\begin{tabular}{l l l l l l l l }
$\mathcal{B}_{2}$:
& \{0, 3, 12\}& \{1, 4, $b$\}& \{5, 6, 11\}& \{8, 14, $a$\}& \{9, 13, 17\}& \{10, 15, 16\}\\
    & \{0, 4, 5\}& \{3, 7, $a$\}& \{6, 12, 17\}& \{8, 10, 13\}& \{9, 16, $b$\}& \{11, 14, 15\}\\
    & \{0, 7, 17\}& \{1, 10, 11\}& \{3, 14, 16\}& \{4, 6, 13\}& \{5, 8, $b$\}& \{9, 12, $a$\}\\
    & \{0, 10, $b$\}& \{1, 5, 14\}& \{3, 8, 17\}& \{4, 11, 12\}& \{6, 7, 16\}& \{13, 15, $a$\}\\
    & \{0, 8, 16\}& \{1, 7, 12\}& \{3, 11, 13\}& \{5, 9, 15\}& \{6, 10, $a$\}& \{14, 17, $b$\}\\
    & \{1, 15, 17\}& \{3, 5, 10\}& \{4, 16, $a$\}& \{6, 9, 14\}& \{7, 8, 11\}& \{12, 13, $b$\}\\
    & \{0, 9, 11\}& \{1, 13, 16\}& \{3, 6, $b$\}& \{4, 7, 14\}& \{5, 17, $a$\}& \{8, 12, 15\}\\
    & \{0, 1, $a$\}& \{4, 8, 9\}& \{7, 15, $b$\}& \{10, 12, 14\}& \{11, 16, 17\}\\
        & \{0, 13, 14\}& \{1, 6, 8\}& \{3, 4, 15\}& \{5, 12, 16\}& \{7, 9, 10\}\\
    & \{0, 6, 15\}& \{1, 3, 9\}& \{4, 10, 17\}& \{5, 7, 13\}\\
\end{tabular}
\end{center}}
We list below the block sets of the derived designs $\mathrm{gcGDD}(2,3,18)$ of type $2^{9}$ at the point $a$ and $b$; the blocks in each row form a $\mathrm{PPC}$ of $I_{18}\setminus G$ for some $G\in\mathcal{G}$. Note that $\mathcal{B}_{a}$ and $\mathcal{B}_{b}$ have the same special group $\{0,9\}$.
{\small\begin{center}
\begin{tabular}{l l l l l l l l }
$\mathcal{B}_{a}$:
    & \{0, 13, 16\}& \{2, 8, 14\}& \{3, 4, 11\}& \{5, 6, 9\}& \{7, 15, 17\}\\
    & \{0, 5, 10\}& \{1, 16, 17\}& \{4, 9, 15\}& \{6, 7, 8\}& \{12, 13, 14\}\\
    & \{0, 7, 11\}& \{1, 4, 6\}& \{2, 5, 17\}& \{8, 9, 13\}& \{10, 14, 16\}\\
    & \{0, 6, 12\}& \{1, 14, 15\}& \{2, 3, 7\}& \{5, 8, 11\}& \{9, 10, 17\}\\
    & \{0, 1, 2\}& \{3, 9, 16\}& \{4, 8, 10\}& \{6, 13, 17\}& \{11, 12, 15\}\\
    & \{0, 3, 17\}& \{1, 5, 12\}& \{2, 4, 16\}& \{7, 9, 14\}& \{10, 11, 13\}\\
    & \{0, 8, 15\}& \{1, 9, 11\}& \{2, 6, 10\}& \{3, 5, 13\}& \{4, 12, 17\}\\
    & \{0, 4, 14\}& \{1, 7, 13\}& \{2, 9, 12\}& \{3, 10, 15\}& \{6, 11, 16\}\\
    & \{2, 13, 15\}& \{3, 6, 14\}& \{4, 5, 7\}& \{8, 12, 16\}\\
    & \{1, 3, 8\}& \{5, 15, 16\}& \{7, 10, 12\}& \{11, 14, 17\}\\
\end{tabular}
\end{center}

\begin{center}
\begin{tabular}{l l l l l l l l }
$\mathcal{B}_{b}$:
    & \{0, 14, 15\}& \{2, 12, 13\}& \{4, 7, 9\}& \{5, 11, 17\}& \{6, 8, 16\}\\ 
    & \{0, 6, 13\}& \{1, 9, 14\}& \{3, 4, 5\}& \{7, 8, 10\}& \{15, 16, 17\}\\ 
    & \{0, 4, 17\}& \{1, 5, 7\}& \{2, 9, 16\}& \{8, 11, 14\}& \{10, 13, 15\}\\
    & \{0, 1, 8\}& \{2, 14, 17\}& \{3, 9, 15\}& \{5, 6, 10\}& \{11, 12, 16\}\\
    & \{0, 7, 12\}& \{1, 13, 17\}& \{2, 3, 6\}& \{4, 8, 15\}& \{9, 10, 11\}\\
    & \{0, 5, 16\}& \{1, 2, 4\}& \{3, 10, 14\}& \{7, 11, 13\}& \{8, 9, 12\}\\
    & \{0, 2, 10\}& \{1, 11, 15\}& \{3, 8, 13\}& \{4, 12, 14\}& \{6, 9, 17\}\\
    & \{0, 3, 11\}& \{1, 6, 12\}& \{2, 7, 15\}& \{4, 10, 16\}& \{5, 9, 13\}\\
    & \{3, 7, 17\}& \{4, 6, 11\}& \{5, 12, 15\}& \{13, 14, 16\}\\
    & \{1, 3, 16\}& \{2, 5, 8\}& \{6, 7, 14\}& \{10, 12, 17\}\\
\end{tabular}
\end{center}}
\qed

{\Lemma\label{Dcqs20_2} There exists a $\mathrm{^{2}cDCQS}(2^{9}:2)$.}
\proof Let $X=I_{18}\cup\{a,b\}$ and $\mathcal{G} = \{\{i,i+9\}: 0\leq i \leq 8\}$.
We construct a $\mathrm{CQS}(2^{9}:2)$ on point set $X$ with stem $S=\{a,b\}$ and group set $\mathcal{G}$.
Let $\mathcal{H}$ be a permutation group generated by the following permutation:
$$ (0,3,6,9,12,15)(1,4,7,10,13,16)(2,5,8,11,14,17)(a)(b).$$
Our $\mathrm{CQS}(2^{9}:2)$ $(X,S,\mathcal{G},\mathcal{B})$ admits the automorphism group $\mathcal{H}$ and the base blocks are listed as follows. (The underlined blocks belong to short orbits of length $2$.)
\begin{center}
\begin{tabular}{l l l l l l l l }
\{0, 1, 2, $a$\}& \{0, 1, 3, 6\}& \{0, 1, 4, 5\}& \{0, 1, 7, 9\}& \{0, 1, 8, $b$\}\\
\{0, 1, 10, 15\}& \{0, 1, 11, 14\}& \{0, 1, 12, 17\}& \{0, 1, 13, 16\}& \{0, 2, 3, 5\}\\
\{0, 2, 4, 8\}& \{0, 2, 6, 11\}& \{0, 2, 7, 16\}& \{0, 2, 9, 12\}& \{0, 2, 10, $b$\}\\
\{0, 2, 13, 14\}& \{0, 3, 7, $a$\}& \{0, 3, 8, 12\}& \{0, 3, 10, 17\}& \{0, 3, 14, $b$\}\\
\{0, 4, 6, $b$\}& \{0, 4, 7, 12\}& \{0, 4, 9, 17\}& \{0, 4, 10, 14\}& \{0, 4, 11, 16\}\\
\{0, 5, 8, 16\}& \{0, 5, 10, $a$\}& \{0, 5, 11, 13\}& \{0, 5, 17, $b$\}& {\underline{\{0, 6, 12, $a$\}}}\\
\{0, 6, 16, 17\}& \{0, 7, 8, 10\}& \{0, 7, 11, 17\}& \{0, 7, 13, $b$\}& \{0, 8, 11, $a$\}\\
\{0, 8, 14, 17\}&  \{0, 13, 17, $a$\}& \{0, 14, 16, $a$\}& \{1, 2, 5, $b$\}&\{1, 2, 7, 10\}\\
\{1, 2, 11, 13\}& \{1, 2, 14, 17\}& \{1, 4, 8, 10\}& \{1, 4, 11, $a$\}&\{1, 4, 17, $b$\}\\
 \{1, 5, 8, 14\}& {\underline{\{1, 7, 13, $a$\}}}&{\underline{\{2, 8, 14, $a$\}}}\\
\end{tabular}
\end{center}%

The resolutions of derived designs at $0,1,2,a$ and $b$ are listed in Appendix C. Thus, we obtain a
$\mathrm{^{2}cDCQS}(2^{9}:2)$.
\qed
{\Lemma \label{RDGDD20} There exists an $\mathrm{RDGDD}(3,4, 10\{2\})$.}

\proof Let $X= I_{20}$ and $\mathcal{G} = \{\{i,i+10\}: 0\leq i \leq 9\}$. We construct a $\mathrm{GDD}(3,4, 20)$ of type $2^{10}$ on $X$ with group set $\mathcal{G}$. Let $\mathcal{H}$ be a permutation group generated by the following permutation:
$$(0,1,2,3,4,5,6,7,8,9)(10,11,12,13,14,15,16,17,18,19).$$
Our $\mathrm{GDD}(3,4,20)$ $(X,\mathcal{G},\mathcal{B})$ admits the automorphism group $\mathcal{H}$ and the base blocks are listed as follows:
\begin{center}
\begin{tabular}{l l l l l l l l }
\{0, 1, 2, 6\}&\{0, 1, 3, 4\}&\{0, 1, 12, 16\}&\{0, 1, 13, 14\}&\{0, 1, 15, 19\}\\
\{0, 1, 17, 18\}&\{0, 2, 4, 17\}&\{0, 2, 5, 18\}&\{0, 2, 7, 14\}&\{0, 2, 11, 16\}\\
\{0, 2, 13, 19\}&\{0, 3, 6, 18\}&\{0, 3, 11, 14\}&\{0, 3, 12, 19\}&\{0, 4, 11, 13\}\\
\{0, 4, 15, 16\}&\{0, 4, 18, 19\}&\{0, 5, 11, 19\}&\{0, 12, 14, 17\}&\{0, 12, 15, 18\}\\
\{0, 13, 15, 17\}&\{0, 13, 16, 18\}&\{10, 11, 12, 16\}&\{10, 11, 13, 14\}\\
\end{tabular}
\end{center}
For any $x\in G$, $G\in \mathcal{G}$, we can check that the derived design at the point $x$ of the $\mathrm{GDD}$ is an $\mathrm{RGDD}(2,3, 18)$ of type $2^{9}$ with group set $\mathcal{G}\setminus \{G\}$. We only need to list the resolutions of the derived designs at $0$ and $10$ to prove; the blocks in each row form a $\mathrm{PC}$ of $X\setminus \{0,10\}$.
{\small\begin{center}
\begin{tabular}{l l l l l l l l }
$\mathcal{B}_{0}$: & \{1, 2, 6\}&\{3, 7, 15\}&\{4, 18, 19\}&\{5, 14, 16\}&\{8, 11, 17\}&\{9, 12, 13\}\\
          &\{1, 3, 4\}&\{2, 7, 14\}&\{5, 8, 12\}&\{6, 17, 19\}&\{9, 11, 15\}&\{13, 16, 18\}\\
          &\{1, 7, 8\}&\{2, 11, 16\}&\{3, 12, 19\}&\{4, 5, 6\}&\{9, 14, 18\}&\{13, 15, 17\}\\
          &\{1, 12, 16\}&\{2, 13, 19\}&\{3, 5, 17\}&\{4, 8, 9\}&\{6, 14, 15\}&\{7, 11, 18\}\\
          &\{1, 5, 9\}&\{2, 4, 17\}&\{3, 11, 14\}&\{6, 8, 13\}&\{7, 16, 19\}&\{12, 15, 18\}\\
          &\{1, 17, 18\}&\{2, 3, 9\}&\{4, 15, 16\}&\{5, 7, 13\}&\{6, 11, 12\}&\{8, 14, 19\}\\
          &\{1, 15, 19\}&\{2, 5, 18\}&\{3, 8, 16\}&\{4, 11, 13\}&\{6, 7, 9\}&\{12, 14, 17\}\\
          &\{1, 13, 14\}&\{2, 8, 15\}&\{3, 6, 18\}&\{4, 7, 12\}&\{5, 11, 19\}&\{9, 16, 17\}\\
\end{tabular}
\end{center}

\begin{center}
\begin{tabular}{l l l l l l l l }
$\mathcal{B}_{10}$: & \{1, 2, 16\}&\{3, 15, 17\}&\{4, 9, 18\}&\{5, 6, 14\}&\{7, 8, 11\}&\{12, 13, 19\}\\
          &\{1, 3, 14\}&\{2, 4, 7\}&\{5, 9, 11\}&\{6, 13, 18\}&\{8, 12, 15\}&\{16, 17, 19\}\\
          &\{1, 4, 13\}&\{2, 5, 8\}&\{3, 16, 18\}&\{6, 9, 17\}&\{7, 12, 14\}&\{11, 15, 19\}\\
          &\{1, 5, 19\}&\{2, 15, 18\}&\{3, 6, 8\}&\{4, 12, 17\}&\{7, 9, 16\}&\{11, 13, 14\}\\
          &\{1, 6, 12\}&\{2, 3, 19\}&\{4, 5, 16\}&\{7, 13, 15\}&\{8, 9, 14\}&\{11, 17, 18\}\\
          &\{1, 8, 17\}&\{2, 9, 13\}&\{3, 5, 7\}&\{4, 6, 15\}&\{11, 12, 16\}&\{14, 18, 19\}\\
          &\{1, 9, 15\}&\{2, 14, 17\}&\{3, 4, 11\}&\{5, 12, 18\}&\{6, 7, 19\}&\{8, 13, 16\}\\
          &\{1, 7, 18\}&\{2, 6, 11\}&\{3, 9, 12\}&\{4, 8, 19\}&\{5, 13, 17\}&\{14, 15, 16\}\\
\end{tabular}
\end{center}}
\qed

We record a useful result on $\mathrm{RDS}$s.

\begin{Lemma}[\label{RDS}{\cite{RDS2,RDS}}]
There exists an $\mathrm{RDS}(3,q+1,q^{n}+1)$ for any prime power $q$ and any positive integer $n$.
\end{Lemma}

{\Theorem\label{result_1} There exists an $\mathrm{mcDSQS}(2\cdot 9^{n}+2)$ for any positive integer $n$.}
\proof For any positive integer $n$, apply Construction \ref{C1} with an $\mathrm{RDS}(3,10,9^{n}+1)$, which exists by Lemma \ref{RDS}. Since a $\mathrm{^{1}cDCQS}(2^{9}:2)$ exists by Lemma \ref{Dcqs20_1}, a $\mathrm{^{2}cDCQS}(2^{9}:2)$ exists by Lemma \ref{Dcqs20_2} and an $\mathrm{RDGDD}(3,4,10\{2\})$ exists by Lemma \ref{RDGDD20}, we can obtain a  $\mathrm{^{1}cDCQS}(2^{9^{n}}: 2)$. An $\mathrm{RDSQS}(4)$ exists obviously. Hence, an $\mathrm{mcDSQS}(2\cdot 9^{n}+2)$ exists by Construction \ref{cqs-wrdsqs1}.
\qed

\begin{Theorem}\label{result_3}
\begin{itemize}
 \item [$(1)$] Let $n\in\{2^{2m+1}+2:m\geq 0\}$. Then $\frac{n}{2}$ is the smallest $q$ such that a diameter perfect $(n,\frac{1}{4}\tbinom{n}{3},6;4)_{q}$ code exists.
 \item [$(2)$] Let $n \in\{ 2\cdot 9^{m}+2:m\geq 1\}\cup\{26,32\}$. Then $\frac{n}{2}+1$ is the smallest $q$ such that a diameter perfect $(n,\frac{1}{4}\tbinom{n}{3},6;4)_{q}$ code exists.
 \end{itemize}
\end{Theorem}
\proof Combine Theorem \ref{result_2} and Corollary \ref{shi-c1} to obtain the desired result in (1). Combine Lemma \ref{MIN_Dsqs}, Theorem \ref{result_1} and Corollary \ref{shi-c1} to obtain the result in (2).
\qed


\section{Summary and open problems}

In this paper, we considered the constructions and existence of $\mathrm{mcDSQS}(v)$s. There has been an amount of  work on RDSQS$(6n+4)$s,
 while the study on $\mathrm{mcDSQS}(6n+2)$s has just begun in this paper. We generalized a number of existing recursive constructions about $\mathrm{RDSQS}$s to produce
  general $\mathrm{mcDSQS}$s. 
In particular, a  construction for $\mathrm{gcDCQS}$s is developed, which is also new even for $\mathrm{RDCQS}$s, see Construction \ref{C4}.
Moreover, we demonstrated several constructions aiming only at $\mathrm{mcDSQS}(6n+2)$s (in Subsection 2.2). As the main results, we obtained a new infinite family of $\mathrm{RDSQS}(6n+4)$s and displayed a new infinite family of $\mathrm{LKTS}$s. (More families of $\mathrm{RDSQS}$s and thus more $\mathrm{LKTS}$s could be produced by substituting known orders to the product construction \cite{Zhou-lkts,yuan_olkts}, but we did not study this in depth to avoid distraction.)
We obtained the first infinite family of $\mathrm{mcDSQS}(6n+2)$s and presented their applications in diameter perfect constant-weight codes. Because the chromatic index of an STS$(v)$ equals $(v+1)/2$ only when $v\equiv 1\pmod 6$ and $v\ge 19$, direct constructions for  auxiliary designs to build $\mathrm{mcDSQS}(6n+2)$s
usually require significant computation.
For future work, we pose the following problems.\vspace{-0.2cm}
\begin{itemize}
\item[(1)] Find more examples of $\mathrm{mcDSQS}(v)$s, for example, an mcDSQS$(38)$ and an RDSQS$(52)$ (the smallest unknown mcDSQS$(6n+2)$s and the smallest unknown
 RDSQS$(6n+4)$s, respectively).\vspace{-0.2cm}
\item[(2)] Find a $\mathrm{gcDCQS}(v^4:0)$ for $v=20,32$. Apply Corollary \ref{Rmul}. If a $\mathrm{gcDCQS}(20^4:0)$ exists, then there exists an $\mathrm{mcDSQS}(5\cdot 4^{n})$ for any integer $n\geq 3$. In the same way, the existence of a $\mathrm{gcDCQS}(32^4:0)$ will yield an $\mathrm{mcDSQS}(2^{2n+1})$ for any  integer $n\geq 4$.\vspace{-0.2cm}
\item[(3)] Construct a $\mathrm{gcDCQS}(18^{4}:2)$ and an $\mathrm{FDGDD}(3,4,5\{18\})$ or  $\mathrm{RDGDD}(3,4,5\{18\})$. By Construction \ref{C3} or \ref{C2} and Lemma \ref{RDS}, the existence of an $\mathrm{mcDSQS}(18\cdot 4^n+2)$ would be obtained by filling with an $\mathrm{mcDSQS}(20)$.\vspace{-0.2cm}
\item[(4)] Develop more effective recursive constructions for $\mathrm{mcDSQS}(v)$s.
\end{itemize}


\appendix
\section{The derived designs in Lemma \ref{CQS}}
{\noindent For the given $\mathrm{CQS}(8^{4}:2)$ in Lemma \ref{CQS}, we list below the block set of the derived $\mathrm{KTS}(33,9)$ at the point $0$.
The first four parts form $\mathrm{PC}$s of $\mathbb{Z}_{32}\setminus G_0$; other parts form $\mathrm{PC}$s of $(\mathbb{Z}_{32}\cup S)\setminus \{0\}$.}
{\small\begin{center}
\begin{tabular}{l l l l l l l l }
$\mathcal{B}_{0}$: &\{1, 7, 9\}& \{2, 13, 22\}& \{3, 17, 18\}& \{5, 11, 26\}& \{6, 21, 27\}& \{10, 19, 30\}\\
\smallskip
& \{14, 15, 29\}& \{23, 25, 31\}\\
&\{1, 10, 25\}& \{2, 3, 19\}& \{5, 9, 14\}& \{6, 11, 17\}& \{7, 22, 31\}& \{13, 29, 30\}\\
\smallskip
& \{15, 21, 26\}& \{18, 23, 27\}\\
&\{1, 14, 19\}& \{2, 11, 27\}& \{3, 9, 26\}& \{5, 21, 30\}& \{6, 23, 29\}& \{7, 10, 15\}\\
\smallskip
& \{13, 18, 31\}& \{17, 22, 25\}\\
&\{1, 15, 18\}& \{2, 25, 26\}& \{3, 10, 27\}& \{5, 22, 29\}& \{6, 7, 30\}& \{9, 19, 21\}\\
\smallskip
& \{11, 13, 23\}& \{14, 17, 31\}\\
&\{1, 31, $b$\}& \{2, 4, 18\}& \{3, 21, 24\}& \{5, 20, 25\}& \{6, 9, 15\}& \{7, 12, 27\}\\
\smallskip
& \{8, 11, 29\}& \{10, 16, 22\}& \{13, 19, $a$\}& \{14, 28, 30\}& \{17, 23, 26\}\\
&\{1, 24, 26\}& \{2, 17, $b$\}& \{3, 11, 14\}& \{4, 6, 10\}& \{5, 23, 28\}& \{7, 13, 21\}\\
\smallskip
& \{8, 27, 30\}& \{9, 20, 22\}& \{12, 15, 19\}& \{16, 29, 31\}& \{18, 25, $a$\}\\
&\{1, 21, 22\}& \{2, 16, 30\}& \{3, 29, $b$\}& \{4, 7, 19\}& \{5, 8, 15\}& \{6, 18, 20\}\\
\smallskip
& \{9, 23, $a$\}& \{10, 11, 31\}& \{12, 14, 26\}& \{13, 25, 28\}& \{17, 24, 27\}\\
&\{1, 17, 30\}& \{2, 15, 31\}& \{3, 8, 25\}& \{4, 14, 22\}& \{5, 27, $b$\}& \{6, 16, 26\}\\
\smallskip
& \{7, 24, 29\}& \{9, 12, 13\}& \{10, 18, 28\}& \{11, 21, $a$\}& \{19, 20, 23\}\\
&\{1, 2, $a$\}& \{3, 5, 13\}& \{4, 17, 29\}& \{6, 19, $b$\}& \{7, 11, 28\}& \{8, 14, 21\}\\
\smallskip
& \{9, 24, 31\}& \{10, 20, 26\}& \{12, 18, 30\}& \{15, 22, 23\}& \{16, 25, 27\}\\
&\{1, 5, 6\}& \{2, 12, 21\}& \{3, 4, 23\}& \{7, 25, $b$\}& \{8, 10, 13\}& \{9, 28, 29\}\\
\smallskip
& \{11, 20, 30\}& \{14, 16, 18\}& \{15, 17, $a$\}& \{19, 22, 24\}& \{26, 27, 31\}\\
&\{1, 11, 12\}& \{2, 10, 29\}& \{3, 7, 20\}& \{4, 5, 31\}& \{6, 14, 25\}& \{8, 19, 26\}\\
\smallskip
& \{9, 18, $b$\}& \{13, 15, 16\}& \{17, 21, 28\}& \{22, 27, $a$\}& \{23, 24, 30\}\\
&\{1, 20, 29\}& \{2, 7, 23\}& \{3, 6, $a$\}& \{4, 26, 30\}& \{5, 12, 17\}& \{8, 18, 22\}\\
\smallskip
& \{9, 11, 16\}& \{10, 21, $b$\}& \{13, 14, 27\}& \{15, 24, 25\}& \{19, 28, 31\}\\
&\{1, 4, 13\}& \{2, 6, 28\}& \{3, 12, 31\}& \{5, 18, 19\}& \{7, 8, 17\}& \{9, 25, 30\}\\
\smallskip
& \{10, 14, 24\}& \{11, 22, $b$\}& \{15, 20, 27\}& \{16, 21, 23\}& \{26, 29, $a$\}\\
&\{1, 8, 23\}& \{2, 14, 20\}& \{3, 15, 28\}& \{4, 21, 25\}& \{5, 7, 16\}& \{6, 12, 22\}\\
\smallskip
& \{9, 10, 17\}& \{11, 18, 24\}& \{13, 26, $b$\}& \{19, 27, 29\}& \{30, 31, $a$\}\\
&\{1, 27, 28\}& \{2, 8, 9\}& \{3, 22, 30\}& \{4, 11, 15\}& \{5, 10, $a$\}& \{6, 13, 24\}\\
\smallskip
& \{7, 18, 26\}& \{12, 25, 29\}& \{14, 23, $b$\}& \{16, 17, 19\}& \{20, 21, 31\}\\
&\{1, 3, 16\}& \{2, 5, 24\}& \{4, 9, 27\}& \{6, 8, 31\}& \{7, 14, $a$\}& \{10, 12, 23\}\\
\smallskip
& \{11, 19, 25\}& \{13, 17, 20\}& \{15, 30, $b$\}& \{18, 21, 29\}& \{22, 26, 28\}\\
\end{tabular}
\end{center}
}
We list below the block set of the derived $\mathrm{KF}(8^{4})$ at $a$; the blocks in each part form a $\mathrm{PC}$ of $\mathbb{Z}_{32}\setminus G$ for some $G\in \mathcal{G}$.
{\small\begin{center}
\begin{tabular}{l l l l l l l l }
$\mathcal{B}_{a}$: &\{1, 3, 18\}& \{2, 17, 19\}& \{5, 7, 22\}& \{6, 21, 23\}\\
\smallskip
& \{9, 11, 26\}& \{10, 25, 27\}& \{13, 15, 30\}& \{14, 29, 31\}\\
&\{1, 10, 19\}& \{2, 11, 25\}& \{3, 17, 26\}& \{5, 14, 23\}\\
\smallskip
& \{6, 15, 29\}& \{7, 21, 30\}& \{9, 18, 27\}& \{13, 22, 31\}\\
&\{1, 11, 22\}& \{2, 13, 23\}& \{3, 14, 25\}& \{5, 15, 26\}\\
\smallskip
& \{6, 17, 27\}& \{7, 18, 29\}& \{9, 19, 30\}& \{10, 21, 31\}\\
&\{1, 14, 27\}& \{2, 15, 21\}& \{3, 9, 22\}& \{5, 18, 31\}\\
\smallskip
& \{6, 19, 25\}& \{7, 13, 26\}& \{10, 23, 29\}& \{11, 17, 30\}\\
&\{0, 3, 6\}& \{2, 7, 12\}& \{4, 26, 31\}& \{8, 11, 14\}\\
\smallskip
& \{10, 15, 20\}& \{16, 19, 22\}& \{18, 23, 28\}& \{24, 27, 30\}\\
&\{0, 7, 14\}& \{2, 3, 4\}& \{6, 24, 31\}& \{8, 15, 22\}\\
\smallskip
& \{10, 11, 12\}& \{16, 23, 30\}& \{18, 19, 20\}& \{26, 27, 28\}\\
&\{0, 22, 27\}& \{2, 28, 31\}& \{3, 8, 30\}& \{4, 7, 10\}\\
\smallskip
& \{6, 11, 16\}& \{12, 15, 18\}& \{14, 19, 24\}& \{20, 23, 26\}\\
&\{0, 30, 31\}& \{2, 20, 27\}& \{3, 10, 28\}& \{4, 11, 18\}\\
\smallskip
& \{6, 7, 8\}& \{12, 19, 26\}& \{14, 15, 16\}& \{22, 23, 24\}\\
&\{0, 9, 23\}& \{1, 15, 24\}& \{3, 5, 20\}& \{4, 19, 21\}\\
\smallskip
& \{7, 16, 25\}& \{8, 17, 31\}& \{11, 13, 28\}& \{12, 27, 29\}\\
&\{0, 11, 21\}& \{1, 7, 20\}& \{3, 13, 24\}& \{4, 17, 23\}\\
\smallskip
& \{5, 16, 27\}& \{8, 19, 29\}& \{9, 15, 28\}& \{12, 25, 31\}\\
&\{0, 13, 19\}& \{1, 12, 23\}& \{3, 16, 29\}& \{4, 15, 25\}\\
\smallskip
& \{5, 11, 24\}& \{7, 17, 28\}& \{8, 21, 27\}& \{9, 20, 31\}\\
&\{0, 15, 17\}& \{1, 16, 31\}& \{3, 12, 21\}& \{4, 13, 27\}\\
\smallskip
& \{5, 19, 28\}& \{7, 9, 24\}& \{8, 23, 25\}& \{11, 20, 29\}\\
&\{0, 1, 2\}& \{4, 5, 6\}& \{8, 9, 10\}& \{12, 13, 14\}\\
\smallskip
& \{16, 17, 18\}& \{20, 21, 22\}& \{24, 25, 26\}& \{28, 29, 30\}\\
&\{0, 5, 10\}& \{1, 4, 30\}& \{2, 24, 29\}& \{6, 9, 12\}\\
\smallskip
& \{8, 13, 18\}& \{14, 17, 20\}& \{16, 21, 26\}& \{22, 25, 28\}\\
&\{0, 18, 25\}& \{1, 8, 26\}& \{2, 9, 16\}& \{4, 22, 29\}\\
\smallskip
& \{5, 12, 30\}& \{6, 13, 20\}& \{10, 17, 24\}& \{14, 21, 28\}\\
&\{0, 26, 29\}& \{1, 6, 28\}& \{2, 5, 8\}& \{4, 9, 14\}\\
\smallskip
& \{10, 13, 16\}& \{12, 17, 22\}& \{18, 21, 24\}& \{20, 25, 30\}\\
\end{tabular}
\end{center}
}

\section{An $\mathrm{mcDSQS}(26)$ and an $\mathrm{mcDSQS}(32)$ in Lemma \ref{MIN_Dsqs}}
{\noindent We first construct a cyclic $\mathrm{SQS}(26)$ on point set $\mathbb{Z}_{26}$. The base blocks are listed as follows. (The underlined blocks belong to short orbits of length $13$.)}
\begin{center}
\begin{tabular}{l l l l l l l l }
\{0, 1, 2, 4\}& \{0, 1, 5, 23\}& \{0, 1, 6, 22\}& \{0, 1, 7, 21\}& \{0, 1, 8, 20\}& \{0, 1, 9, 15\}\\
\{0, 1, 10, 12\}& \{0, 1, 11, 19\}& \underline{\{0, 1, 13, 14\}}& \{0, 1, 16, 18\}& \{0, 1, 17, 24\}& \{0, 2, 5, 7\}\\
\{0, 2, 6, 9\}& \{0, 2, 8, 13\}& \{0, 2, 12, 22\}& \{0, 2, 14, 18\}& \{0, 2, 15, 20\}& \{0, 3, 6, 10\}\\
\{0, 3, 9, 18\}& \{0, 3, 11, 15\}& \{0, 3, 12, 17\}& \underline{\{0, 3, 13, 16\}}& \{0, 3, 14, 21\}& \{0, 4, 9, 19\}\\
\{0, 4, 11, 21\}& \underline{\{0, 4, 13, 17\}}& \underline{\{0, 6, 13, 19\}}\\
\end{tabular}
\end{center}
The above cyclic $\mathrm{SQS}(26)$ is an $\mathrm{mcDSQS}(26)$.
We list below the block set of the derived design at the point $0$; the blocks in each part form a $\mathrm{PPC}$ of $\mathbb{Z}_{26}\setminus\{0\}$.
{\small\begin{center}
\begin{tabular}{l l l l l l l l l l}
\{1, 2, 4\}&\{3, 5, 24\}&\{6, 7, 14\}&\{8, 10, 22\}\\
\smallskip
\{9, 17, 20\}&\{12, 18, 19\}&\{13, 15, 21\}&\{16, 23, 25\}\\
\{1, 3, 25\}&\{2, 5, 7\}&\{4, 13, 17\}&\{6, 15, 23\}\\
\smallskip
\{8, 9, 24\}&\{10, 14, 16\}&\{11, 12, 20\}&\{18, 21, 22\}\\
\{1, 5, 23\}&\{2, 3, 19\}&\{4, 6, 16\}&\{7, 8, 18\}\\
\smallskip
\{9, 12, 21\}&\{10, 20, 24\}&\{11, 14, 22\}&\{15, 17, 25\}\\
\{1, 6, 22\}&\{2, 10, 11\}&\{3, 12, 17\}&\{4, 15, 18\}\\
\smallskip
\{5, 9, 16\}&\{7, 13, 20\}&\{8, 14, 25\}&\{19, 21, 24\}\\
\{1, 7, 21\}&\{2, 14, 18\}&\{3, 20, 22\}&\{4, 5, 10\}\\
\smallskip
\{6, 11, 24\}&\{8, 15, 16\}&\{12, 13, 25\}&\{17, 19, 23\}\\
\{1, 8, 20\}&\{2, 16, 17\}&\{3, 11, 15\}&\{4, 12, 14\}\\
\smallskip
\{5, 21, 25\}&\{6, 13, 19\}&\{7, 9, 10\}&\{22, 23, 24\}\\
\{1, 9, 15\}&\{2, 21, 23\}&\{3, 13, 16\}&\{4, 7, 24\}\\
\smallskip
\{5, 6, 12\}&\{8, 11, 17\}&\{10, 18, 25\}&\{14, 19, 20\}\\
\{1, 10, 12\}&\{2, 24, 25\}&\{3, 14, 21\}&\{4, 20, 23\}\\
\smallskip
\{5, 8, 19\}&\{6, 17, 18\}&\{7, 11, 16\}&\{9, 13, 22\}\\
\{1, 11, 19\}&\{2, 15, 20\}&\{3, 9, 18\}&\{4, 22, 25\}\\
\smallskip
\{5, 14, 17\}&\{6, 8, 21\}&\{10, 13, 23\}&\{12, 16, 24\}\\
\{1, 16, 18\}&\{2, 12, 22\}&\{3, 7, 23\}&\{4, 9, 19\}\\
\smallskip
\{5, 11, 13\}&\{6, 20, 25\}&\{10, 17, 21\}&\{14, 15, 24\}\\
\{1, 17, 24\}&\{2, 8, 13\}&\{3, 6, 10\}&\{4, 11, 21\}\\
\smallskip
\{5, 18, 20\}&\{7, 12, 15\}&\{9, 14, 23\}&\{16, 19, 22\}\\
\{1, 13, 14\}&\{2, 6, 9\}&\{3, 4, 8\}&\{5, 15, 22\}\\
\smallskip
\{7, 19, 25\}&\{11, 18, 23\}\\
\{7, 17, 22\}&\{8, 12, 23\}&\{9, 11, 25\}&\{10, 15, 19\}\\
\smallskip
\{13, 18, 24\}&\{16, 20, 21\}\\
\end{tabular}
\end{center}
}

Now, we construct a cyclic $\mathrm{SQS}(32)$ on point set $\mathbb{Z}_{32}$. The base blocks are listed as follows. (The underlined blocks belong to short orbits of length $16$ while the double-underlined blocks lie in short orbits of length $8$.)
\begin{center}
\begin{tabular}{l l l l l l l l }
\{0, 1, 2, 25\}& \{0, 1, 3, 29\}& \{0, 1, 4, 8\}& \{0, 1, 5, 27\}& \{0, 1, 6, 9\}& \{0, 1, 7, 30\}\\
\{0, 1, 10, 13\}& \{0, 1, 11, 15\}& \{0, 1, 12, 21\}& \{0, 1, 14, 17\}& \{0, 1, 16, 22\}& \{0, 1, 18, 26\}\\
\{0, 1, 19, 23\}& \{0, 1, 20, 28\}& \{0, 2, 4, 15\}& \{0, 2, 5, 26\}& \{0, 2, 6, 14\}& \{0, 2, 7, 23\}\\
\{0, 2, 8, 29\}& \{0, 2, 10, 21\}& \{0, 2, 11, 19\}& \{0, 2, 12, 20\}& \underline{\{0, 2, 16, 18\}}& \{0, 2, 17, 22\}\\
\{0, 2, 24, 27\}& \{0, 3, 12, 17\}& \{0, 3, 13, 25\}& \{0, 3, 14, 20\}& \{0, 3, 15, 28\}& \{0, 3, 16, 21\}\\
\{0, 4, 9, 18\}& \{0, 4, 10, 16\}& \{0, 4, 11, 25\}& \{0, 4, 17, 23\}& \{0, 4, 20, 27\}& \{0, 5, 12, 18\}\\
\{0, 5, 13, 22\}& \{0, 6, 13, 23\}& \{0, 7, 14, 22\}&\underline{\underline{\{0, 8, 16, 24\}}}
\end{tabular}
\end{center}
The above cyclic $\mathrm{SQS}(32)$ is an $\mathrm{mcDSQS}(32)$. We list below the block set of the derived design at the point $0$; the blocks in each part form a $\mathrm{PPC}$ of $\mathbb{Z}_{32}\setminus\{0\}$.
{\small\begin{center}
\begin{tabular}{l l l l l l l l l l}
\{2, 3, 9\}& \{4, 7, 19\}& \{5, 12, 18\}& \{6, 16, 17\}& \{8, 20, 22\}\\\smallskip
\{10, 14, 31\}& \{11, 28, 30\}& \{13, 15, 24\}& \{21, 27, 29\}& \{23, 25, 26\}\\

\{1, 3, 29\}& \{4, 5, 24\}& \{6, 12, 28\}& \{7, 11, 18\}& \{8, 17, 27\}\\\smallskip
 \{9, 20, 21\}& \{10, 15, 23\}& \{13, 16, 31\}& \{14, 19, 26\}& \{22, 25, 30\}\\

\{1, 7, 30\}& \{2, 4, 15\}& \{5, 14, 28\}& \{6, 29, 31\}& \{8, 18, 25\}\\\smallskip
 \{9, 23, 27\}& \{10, 17, 24\}& \{11, 13, 21\}& \{12, 19, 22\}& \{16, 20, 26\}\\

\{1, 2, 25\}& \{3, 5, 11\}& \{6, 18, 21\}& \{7, 14, 22\}& \{8, 26, 28\}\\\smallskip
 \{9, 19, 24\}& \{10, 12, 27\}& \{13, 17, 20\}& \{15, 16, 29\}& \{23, 30, 31\}\\

\{1, 5, 27\}& \{2, 7, 23\}& \{3, 14, 20\}& \{4, 10, 16\}& \{8, 19, 30\}\\\smallskip
\{9, 13, 26\}& \{11, 22, 24\}& \{12, 25, 29\}& \{15, 21, 31\}& \{17, 18, 28\}\\

\{1, 19, 23\}& \{2, 28, 31\}& \{3, 16, 21\}& \{4, 12, 30\}& \{5, 9, 25\}\\\smallskip
\{6, 22, 26\}& \{7, 13, 27\}& \{8, 14, 15\}& \{11, 17, 29\}& \{18, 20, 24\}\\

\{1, 10, 13\}& \{2, 5, 26\}& \{3, 4, 6\}& \{7, 12, 16\}& \{8, 21, 23\}\\\smallskip
\{9, 17, 30\}& \{11, 14, 27\}& \{18, 22, 31\}& \{19, 20, 29\}& \{24, 25, 28\}\\

\{1, 18, 26\}& \{2, 6, 14\}& \{3, 13, 25\}& \{4, 21, 22\}& \{5, 7, 29\}\\
\smallskip
\{8, 16, 24\}& \{9, 10, 28\}& \{11, 12, 23\}& \{15, 20, 30\}& \{19, 27, 31\}\\

\{1, 14, 17\}& \{2, 10, 21\}& \{3, 18, 19\}& \{4, 20, 27\}& \{5, 13, 22\}\\\smallskip
\{6, 8, 11\}& \{7, 15, 25\}& \{9, 12, 31\}& \{16, 23, 28\}& \{26, 29, 30\}\\

\{1, 4, 8\}& \{2, 12, 20\}& \{3, 7, 31\}& \{5, 15, 17\}& \{6, 27, 30\}\\\smallskip
\{9, 11, 16\}& \{10, 22, 29\}& \{13, 19, 28\}& \{14, 18, 23\}& \{21, 24, 26\}\\

\{1, 16, 22\}& \{2, 8, 29\}& \{3, 15, 28\}& \{4, 9, 18\}& \{5, 21, 30\}\\\smallskip
\{6, 13, 23\}& \{7, 17, 26\}& \{10, 19, 25\}& \{11, 20, 31\}& \{12, 14, 24\}\\

\{1, 6, 9\}& \{2, 24, 27\}& \{3, 22, 23\}& \{4, 11, 25\}& \{5, 8, 31\}\\\smallskip
\{7, 10, 20\}& \{12, 15, 26\}& \{13, 18, 29\}& \{14, 16, 30\}& \{17, 19, 21\}\\

\{1, 12, 21\}& \{2, 16, 18\}& \{3, 24, 30\}& \{4, 13, 14\}& \{5, 20, 23\}\\\smallskip
\{6, 15, 19\}& \{7, 8, 9\}& \{10, 11, 26\}& \{17, 25, 31\}& \{22, 27, 28\}\\

\{1, 24, 31\}& \{2, 13, 30\}& \{3, 8, 10\}& \{4, 17, 23\}& \{5, 16, 19\}\\\smallskip
\{6, 20, 25\}& \{7, 21, 28\}& \{9, 14, 29\}& \{15, 18, 27\}\\

\{1, 11, 15\}& \{2, 17, 22\}& \{3, 26, 27\}& \{4, 28, 29\}& \{6, 7, 24\}\\\smallskip
\{8, 12, 13\}& \{10, 18, 30\}& \{14, 21, 25\}\\

\{1, 20, 28\}& \{2, 11, 19\}& \{3, 12, 17\}& \{4, 26, 31\}& \{5, 6, 10\}\\\smallskip
\{9, 15, 22\}& \{16, 25, 27\}& \{23, 24, 29\}\\
\end{tabular}
\end{center}
}

\section{The derived designs in Lemma \ref{Dcqs20_2}}
{\noindent For the given $\mathrm{CQS}(2^{9}:2)$ in Lemma \ref{Dcqs20_2}, we list below the block set of the derived design at the point $x\in\{0,1,2\}$; the blocks in each row form a $\mathrm{PPC}$ of $X\setminus\{x\}$ and the last two rows form $\mathrm{PPC}$s of $I_{18}\setminus G$, where $x\in G$, $G\in \mathcal{G}$.}
{\small\begin{center}
\begin{tabular}{l l l l l l l l }
$\mathcal{B}_{0}$: & \{2, 9, 12\}& \{3, 14, $b$\}& \{4, 15, $a$\}& \{5, 11, 13\}& \{6, 16, 17\}& \{7, 8, 10\}\\
& \{1, 4, 5\}& \{2, 7, 16\}& \{3, 8, 12\}& \{6, 9, 14\}& \{11, 15, $b$\}& \{13, 17, $a$\}\\
& \{1, 12, 17\}& \{2, 3, 5\}& \{4, 6, $b$\}& \{7, 14, 15\}& \{8, 11, $a$\}& \{9, 10, 16\}\\
& \{1, 8, $b$\}& \{2, 15, 17\}& \{3, 9, 11\}& \{4, 7, 12\}& \{6, 10, 13\}& \{14, 16, $a$\}\\
& \{1, 7, 9\}& \{2, 6, 11\}& \{3, 4, 13\}& \{5, 10, $a$\}& \{8, 14, 17\}& \{12, 16, $b$\}\\
& \{1, 13, 16\}& \{2, 10, $b$\}& \{3, 7, $a$\}& \{4, 9, 17\}& \{5, 12, 14\}& \{6, 8, 15\}\\
& \{1, 11, 14\}& \{2, 4, 8\}& \{3, 10, 17\}& \{5, 9, 15\}& \{6, 12, $a$\}& \{7, 13, $b$\}\\
& \{1, 2, $a$\}& \{3, 15, 16\}& \{5, 17, $b$\}& \{8, 9, 13\}& \{10, 11, 12\}\\
& \{1, 3, 6\}& \{4, 10, 14\}& \{5, 8, 16\}& \{7, 11, 17\}& \{12, 13, 15\}\\
& \{1, 10, 15\}& \{2, 13, 14\}& \{4, 11, 16\}& \{5, 6, 7\}\\
\end{tabular}
\end{center}}

{\small\begin{center}
\begin{tabular}{l l l l l l l l }
$\mathcal{B}_{1}$: & \{0, 2, $a$\}& \{3, 5, 10\}& \{4, 13, 14\}& \{6, 12, 16\}& \{7, 8, 17\}& \{9, 11, $b$\}\\
& \{0, 3, 6\}& \{2, 5, $b$\}& \{4, 11, $a$\}& \{8, 13, 15\}& \{9, 16, 17\}& \{10, 12, 14\}\\
& \{0, 4, 5\}& \{2, 3, 9\}& \{6, 14, 15\}& \{7, 12, $b$\}& \{8, 16, $a$\}& \{10, 13, 17\}\\
& \{0, 7, 9\}& \{2, 15, 16\}& \{3, 17, $a$\}& \{4, 8, 10\}& \{5, 11, 12\}& \{6, 13, $b$\}\\
& \{0, 10, 15\}& \{2, 11, 13\}& \{3, 7, 14\}& \{4, 17, $b$\}& \{5, 6, $a$\}& \{8, 9, 12\}\\
& \{0, 8, $b$\}& \{2, 14, 17\}& \{4, 6, 7\}& \{5, 9, 13\}& \{10, 11, 16\}& \{12, 15, $a$\}\\
& \{2, 4, 12\}& \{3, 8, 11\}& \{5, 15, 17\}& \{6, 9, 10\}& \{7, 13, $a$\}& \{14, 16, $b$\}\\
& \{0, 13, 16\}& \{2, 7, 10\}& \{3, 15, $b$\}& \{6, 11, 17\}& \{9, 14, $a$\}\\
& \{0, 11, 14\}& \{2, 6, 8\}& \{3, 12, 13\}& \{4, 9, 15\}& \{5, 7, 16\}\\
& \{0, 12, 17\}& \{3, 4, 16\}& \{5, 8, 14\}& \{7, 11, 15\}\\
\end{tabular}
\end{center}
\begin{center}
\begin{tabular}{l l l l l l l l }
$\mathcal{B}_{2}$:
& \{0, 1, $a$\}& \{3, 11, 17\}& \{4, 9, 14\}& \{5, 7, 8\}& \{6, 10, 16\}& \{12, 13, $b$\}\\
& \{0, 3, 5\}& \{4, 10, 11\}& \{6, 12, 15\}& \{7, 9, 13\}& \{8, 14, $a$\}& \{16, 17, $b$\}\\
& \{0, 6, 11\}& \{1, 5, $b$\}& \{3, 7, 12\}& \{4, 13, 16\}& \{8, 10, 15\}& \{9, 17, $a$\}\\
& \{0, 13, 14\}& \{1, 6, 8\}& \{3, 16, $a$\}& \{4, 7, $b$\}& \{9, 11, 15\}& \{10, 12, 17\}\\
& \{0, 7, 16\}& \{1, 3, 9\}& \{4, 5, 17\}& \{8, 11, 12\}& \{10, 13, $a$\}& \{14, 15, $b$\}\\
& \{0, 15, 17\}& \{1, 11, 13\}& \{3, 8, $b$\}& \{4, 6, $a$\}& \{5, 9, 10\}& \{12, 14, 16\}\\
& \{0, 4, 8\}& \{1, 14, 17\}& \{5, 11, 16\}& \{6, 9, $b$\}& \{7, 15, $a$\}\\
& \{0, 10, $b$\}& \{1, 15, 16\}& \{3, 6, 13\}& \{5, 12, $a$\}& \{7, 11, 14\}\\
& \{1, 4, 12\}& \{3, 10, 14\}& \{5, 13, 15\}& \{6, 7, 17\}& \{8, 9, 16\}\\
& \{0, 9, 12\}& \{1, 7, 10\}& \{3, 4, 15\}& \{5, 6, 14\}& \{8, 13, 17\}\\
\end{tabular}
\end{center}}
We list below the block sets of the derived designs at $a$ and $b$; the blocks in each row form a $\mathrm{PC}$ of $I_{18}$.
{\small\begin{center}
\begin{tabular}{l l l l l l l l }
$\mathcal{B}_{a}$:& \{0, 1, 2\}& \{3, 4, 5\}& \{6, 7, 8\}& \{9, 10, 11\}& \{12, 13, 14\}& \{15, 16, 17\}\\
& \{0, 3, 7\}& \{1, 9, 14\}& \{2, 10, 13\}& \{4, 12, 17\}& \{5, 8, 15\}& \{6, 11, 16\}\\
& \{0, 4, 15\}& \{1, 8, 16\}& \{2, 5, 12\}& \{3, 11, 14\}& \{6, 9, 13\}& \{7, 10, 17\}\\
& \{0, 5, 10\}& \{1, 4, 11\}& \{2, 7, 15\}& \{3, 8, 13\}& \{6, 14, 17\}& \{9, 12, 16\}\\
& \{0, 6, 12\}& \{1, 3, 17\}& \{2, 8, 14\}& \{4, 10, 16\}& \{5, 7, 9\}& \{11, 13, 15\}\\
& \{0, 8, 11\}& \{1, 12, 15\}& \{2, 9, 17\}& \{3, 6, 10\}& \{4, 7, 14\}& \{5, 13, 16\}\\
& \{0, 13, 17\}& \{1, 5, 6\}& \{2, 3, 16\}& \{4, 8, 9\}& \{7, 11, 12\}& \{10, 14, 15\}\\
& \{0, 14, 16\}& \{1, 7, 13\}& \{2, 4, 6\}& \{3, 9, 15\}& \{5, 11, 17\}& \{8, 10, 12\}\\
\end{tabular}
\end{center}
\begin{center}
\begin{tabular}{l l l l l l l l }
$\mathcal{B}_{b}$:& \{0, 1, 8\}& \{2, 12, 13\}& \{3, 4, 11\}& \{5, 15, 16\}& \{6, 7, 14\}& \{9, 10, 17\}\\
& \{0, 2, 10\}& \{1, 9, 11\}& \{3, 5, 13\}& \{4, 12, 14\}& \{6, 8, 16\}& \{7, 15, 17\}\\
& \{0, 3, 14\}& \{1, 6, 13\}& \{2, 16, 17\}& \{4, 10, 15\}& \{5, 9, 12\}& \{7, 8, 11\}\\
& \{0, 4, 6\}& \{1, 14, 16\}& \{2, 3, 8\}& \{5, 7, 10\}& \{9, 13, 15\}& \{11, 12, 17\}\\
& \{0, 5, 17\}& \{1, 3, 15\}& \{2, 4, 7\}& \{6, 10, 12\}& \{8, 9, 14\}& \{11, 13, 16\}\\
& \{0, 7, 13\}& \{1, 2, 5\}& \{3, 6, 17\}& \{4, 9, 16\}& \{8, 12, 15\}& \{10, 11, 14\}\\
& \{0, 11, 15\}& \{1, 7, 12\}& \{2, 6, 9\}& \{3, 10, 16\}& \{4, 5, 8\}& \{13, 14, 17\}\\
& \{0, 12, 16\}& \{1, 4, 17\}& \{2, 14, 15\}& \{3, 7, 9\}& \{5, 6, 11\}& \{8, 10, 13\}\\
\end{tabular}
\end{center}}
\clearpage

\end{document}